\documentclass[12pt]{article}

\usepackage{amsmath,amssymb,latexsym, amsbsy, amsfonts, amscd, amsthm, mathrsfs, xy,comment, xcolor}
\usepackage{multirow}
\usepackage{hyperref}
\usepackage{enumitem, empheq, url}
\usepackage{tikz-cd}
\usepackage{rotating}
\usepackage[bibtex-style]{amsrefs}
\usepackage{ltablex}
\usepackage{array}
\usepackage{mathtools}

\DeclareMathAlphabet\euscr{U}{eus}{m}{n}
\xyoption{all}

\newcommand{\longhookrightarrow}{\lhook\joinrel\longrightarrow}
\newcommand{\longtwoheadrightarrow}{\relbar\joinrel\twoheadrightarrow}

\theoremstyle{plain}

\newtheorem{theorem}{Theorem}[section]
\newtheorem{conjecture}[theorem]{Conjecture}

\newtheorem{lemma}[theorem]{Lemma}

\theoremstyle{definition}
\newtheorem{definition}[theorem]{Definition}

\long\def\symbolfootnote[#1]#2{\begingroup
\def\thefootnote{\fnsymbol{footnote}}\footnote[#1]{#2}\endgroup}

\DeclareMathOperator{\Ann}{Ann}

\def\alg{{\mathrm{alg}}}
\def\an{{\mathrm{an}}}

\def\SL{{\bf SL}}
\def\GL{{\bf GL}}

\def\A{\mathbf{A}}
\def\I{\mathbf{I}}
\def\cA{\mathcal{A}}

\def\N{\mathrm{N}}
\def\1{\mf{1}}

\DeclareMathOperator{\Sku}{SKu}
\DeclareMathOperator{\Sel}{Sel}
\DeclareMathOperator{\Fitt}{Fitt}
\DeclareMathOperator{\Aug}{Aug}

\DeclareMathOperator{\Frac}{Frac}

\DeclareMathOperator{\cyc}{cyc}

 \DeclareMathOperator{\Norm}{Norm}
\DeclareMathOperator{\rec}{rec} 
\DeclareMathOperator{\Hom}{Hom} 

\DeclareMathOperator{\sign}{sign}

\DeclareMathOperator{\ord}{ord}

\DeclareMathOperator{\BS}{BS} \DeclareMathOperator{\Cl}{Cl}

 \DeclareMathOperator{\Frob}{Frob}

\newcommand{\bchi}{\pmb{\chi}}

\newcommand{\mint}{\times\!\!\!\!\!\!\!\int}

\newcommand{\mf}{\mathfrak }

\newcommand{\mscr}{\mathscr }

\def\fa{\mathfrak{a}}
\def\fn{\mathfrak{n}}
\def\fp{\mathfrak{p}}
\def\fq{\mathfrak{q}}

\def\fg{\mathfrak{g}}
\def\fm{\mathfrak{m}}

\def\fl{\mathfrak{l}}
\def\fP{\mathfrak{P}}

\def\fm{\mathfrak{m}}

\def\T{\mathbf{T}}
\def\Z{\mathbf{Z}}

\def\Q{\mathbf{Q}}
\def\C{\mathbf{C}}

\def\R{\mathbf{R}}

\def\bO{\mathbf{O}}

\def\bdf{\begin{defn}}
\def\edf{\end{defn}}

\def\fb{\mathfrak{b}}

\def\Gal{{\rm Gal}}

\def\ab{{\rm ab}}

\def\ram{\text{ram}}
\def\ab{\text{ab}}

\def\sL{{\mscr L}}

\begin{document}
\baselineskip 16.5pt

\title{On the Brumer--Stark Conjecture and Refinements}
\author{Samit Dasgupta \\ Mahesh Kakde}

\maketitle

\begin{abstract}
We state the Brumer--Stark conjecture and motivate it from two perspectives.
Stark's perspective arose in his attempts to generalize the classical Dirichlet class number formula for the leading term of the Dedekind zeta function at $s=1$ (equivalently, $s=0$).  Brumer's perspective arose by generalizing Stickelberger's work regarding the factorization of Gauss sums and the annihilation of class groups of cyclotomic fields.  These viewpoints were synthesized by Tate, who stated the Brumer--Stark conjecture in its current form.

\smallskip
The conjecture considers a totally real field $F$ and a finite abelian CM extension $H/F$.  It states the existence of $p$-units in $H$ whose valuations at places above $p$ are related to the special values of the $L$-functions of the extension $H/F$ at $s=0$.  Essentially equivalently, the conjecture states that a Stickelberger element associated to $H/F$ annihilates the (appropriately smoothed) class group of $H$.  

\smallskip
This conjecture has been refined by many authors in multiple directions.  Notably, Kurihara conjectured that the Stickelberger element lies  in the Fitting ideal of the Pontryagin dual of the class group, and furthermore conjectured an exact formula for this Fitting ideal.  Burns constructed a Selmer group whose Fitting ideal he conjectured to be generated by the Stickelberger element.  Atsuta and Kataoka  conjectured a formula for the Fitting ideal of the class group, rather than its dual.  Rubin stated a higher rank generalization of the Brumer--Stark conjecture.  The first author and his collaborators stated an exact $p$-adic analytic formula for Brumer--Stark units, generalizing conjectures of Gross that give the image of the units under the Artin reciprocity map.

\smallskip
We conclude by stating our results toward these various conjectures and summarizing the proofs.  In particular, we prove the Brumer--Stark conjecture, Rubin's higher rank version, and Kurihara's conjecture, all ``away from 2.'' 
We also prove strong partial results toward Gross's conjecture and the exact $p$-adic analytic formula for Brumer--Stark units.
 The key technique involved in the proofs is Ribet's method.  We demonstrate congruences between Hilbert modular Eisenstein series and cusp forms, and use the associated Galois representations to construct Galois cohomology classes.  These cohomology classes are interpreted in terms of Ritter--Weiss modules, from which results on class groups may be deduced.

\end{abstract}

\maketitle

\section{Background and Motivation}

 Dirichlet's class number formula, conjectured for quadratic fields by Jacobi in 1832 and proven by Dirichlet in 1839, is one of the earliest examples of a relationship between leading terms of $L$-functions and global arithmetic invariants. Let $F$ be a number field with  ring of integers $O_F$. The Dedekind zeta function associated with $F$ is defined as
\[
\zeta_F(s) = \sum_{0 \neq\mathfrak{a} \subset O_F} \N\fa^{-s}, \qquad \qquad \text{Re}(s) > 1,
\]
where $\fa$ runs through the non-zero ideals in $O_F$. The function $\zeta_F(s)$ generalizes Riemann's zeta function and has a meromorphic continuation to the whole complex plane with only a simple pole at $s=1$. Dirichlet's class number formula, which is proved using a ``geometry of numbers'' approach, evaluates the residue at $s=1$:
\[
\lim_{s \rightarrow 1} (s-1)\zeta_F(s) = \frac{2^{r_1} (2 \pi)^{r_2} R_F h_F}{w_F \sqrt{|D_F|}}.
\]
Here  $r_1$ is the number of real embeddings of $F$ and $2r_2$ is the number of complex embeddings of $F$. Further, $h_F$ and $R_F$ denote the class number and regulator (defined below) of $F$, respectively, while $w_F$ denotes the number of roots of unity in $F$ and $D_F$ is the discriminant of $F/\Q$.
The meromorphic function $\zeta_F(s)$ satisfies a functional equation relating $\zeta_F(s)$ and $\zeta_F(1-s)$. Using this functional equation, Dirichlet's class number formula can be restated as giving the leading term of the 
 Taylor expansion of $\zeta_F(s)$ at $s=0$:
\begin{equation} \label{e:dcnf}
\zeta_F(s) =  -\frac{h_F R_F}{w_F}s^{r_1+r_2-1} + O(s^{r_1+r_2}).
\end{equation}

Artin described a theory of $L$-functions generalizing the Dedekind zeta function. Let $G_F$ be the absolute Galois group of $F$. A Dirichlet character for $F$ (or a degree 1 Artin character of $F$) is a homomorphism $\chi\colon G_F \longrightarrow {\C}^{\times}$ with finite image. Class field theory identifies $\chi$ with a function, again denoted by $\chi$, from the set of non-zero ideals of $O_F$ to ${\C}$. Define
\[
L(\chi, s) = \sum_{0 \neq \mathfrak{a} \subset O_F} \chi(\fa) \N\fa^{-s}, \qquad \qquad \text{Re}(s) > 1.
\]
Again $L(\chi,s)$ has a meromorphic continuation to the whole complex plane with only a simple pole at $s=1$ if $\chi$ is trivial. 
If $H/F$ is a Galois extension with finite abelian Galois group $G= \text{Gal}(H/F) $, then we can view any character $\chi \in \hat{G} = \Hom(G, \C^*)$ as a Dirichlet character for $F$, and we have the Artin  factorization formula
\begin{equation}\label{e:factor}
\zeta_H(s) = \prod_{\chi \in \hat{G}} L(\chi, s).
\end{equation}

Dirichlet's class number formula (\ref{e:dcnf}) for the field $H$ gives the leading term of the left hand side of (\ref{e:factor}) at $s=0$.
Stark asked for an analogous formula for $L(\chi, s)$ at $s=0$ for each character $\chi$, thereby giving a canonical factorization of the term $h_HR_H/w_H$.
This led to the formulation of the abelian Stark conjecture, which we state in 
 \S\ref{s:stark}. This statement involves the choice of  places of $F$ that split completely in $H$.   After stating Stark's conjecture,  we restrict in the remainder of the paper to the case that the splitting places of $F$ are finite.  Since the associated $L$-values here are algebraic, one can make progress on the conjectures through $p$-adic techniques such as $p$-adic Galois cohomology.
 To obtain nonzero $L$-values (and hence have nontrivial statements), parity conditions force us to restrict to  the setting that $F$ is a totally real field and $H$ is a CM field.

 Stark's conjecture at finite places has a natural restatement in terms of annihilators of class groups as formulated in the Brumer--Stark conjecture.  We recall the statement and its refinements in \S\ref{s:bs}--\S\ref{s:refine}.  The rest of the paper is taken up in describing the statement and proofs of our results toward the Brumer--Stark conjecture and its refinements.

\section{Stark's conjecture} \label{s:stark} Let us first reformulate Dirichlet's class number formula. 

 For any place $w$ of $F$  we normalize the absolute value $| \ |_w : F_w^* \rightarrow \R$ by 
\[
|u|_w = \left\{ \begin{array}{ll} |u| & \text{ if } w \text{ is real} \\ |u|^2 & \text{ if } w \text{ is complex} \\ \N w^{-\ord_w(u)} & \text{ if } w \text{ is nonarchimedean.}  \end{array} \right. 
\]
For a finite set of places $S$ of $F$, let $X_S$ denote the degree zero subgroup of the free abelian group on $S$.
Let $u_1, \ldots, u_{r_1+r_2-1}$ be a set of generators of the free abelian group $O_F^*/\mu_F$. Let $S_\infty$ be the set of archimedean places of $F$.

The Dirichlet regulator map
\[
O_F^*/\mu_F \rightarrow \R X_{S_\infty}, \qquad
u \mapsto \sum_{w \in S_{\infty}} \log|u|_w \cdot w
\]
induces an isomorphism $\R O_F^* \rightarrow \R X_{S_{\infty}}$. 
Here and throughout, $\R X_{S_\infty}$ denotes  $\R \otimes_{\Z} X_{S_\infty}$, etc.
Let $w_1, \dotsc, w_{r_1+r_2}$ denote the archimedean places of $F$. Then \begin{equation} \label{e:wbasis}
 \{w_i-w_1:  2 \leq i \leq r+1\}, \quad r = r_1+r_2-1\end{equation} is an integral basis of $X_{S_\infty}$. Let $R_F$ be the absolute value of the determinant of the isomorphism between $\R O_F^*$ and $\R X_{S_\infty}$ with respect to the bases $\{ u_1, \ldots , u_{r_1+r_2-1} \}$ and (\ref{e:wbasis}), respectively. Up to a sign, Dirichlet's class number formula can be restated as follows:
\begin{enumerate}
\item[(i)] The rational structure $\Q O_F^*$ on the left hand side corresponds to the structure $\zeta_F^{(r)}(0) \Q X_{S_\infty}$ on the right hand side. 
\item[(ii)] The integral structure $O_F^*/\mu_F$ on the left hand side corresponds to the structure $\zeta_F^{(r)}(0) X_{S_{\infty}}$ on the right hand side.
\end{enumerate}

Motivated by this reformulation we present Stark's conjecture and its  integral refinement due to Rubin. For details see \cite{rubin}. Let $F$ be a number field of degree $n$ and let $H/F$ be a finite Galois extension with $G = \Gal(H/F)$ abelian.  Let $S, T$ be two finite disjoint sets of places of $F$ satisfying the following conditions:
\begin{enumerate}  
\item $S$ contains the sets $S_\infty$ of archimedean places and $S_\ram$ of places ramified in $H$.
\item \label{i:conditions} $T$ contains at least two primes of different residue characteristic or at least one prime of residue characteristic larger than $n+1$, where $n=[F:\Q]$.
\end{enumerate}
For any character $\chi \in \hat{G} = \Hom(G, \C^*),$ define the $S$-depleted, $T$-smoothed $L$-function
\[ L_{S, T}(\chi, s) = \prod_{\fp \not \in S} \frac{1}{1 - \chi(\fp)\N\fp^{-s}} \prod_{\fp \in T}(1 - \chi(\fp) \N\fp^{1-s}), \quad \Re(s) > 1. \]
The function $L_{S,T}(\chi, s)$ extends by analytic continuation to a holomorphic function on $\C$.
 The {\em Stickelberger element} associated to this data is the unique group-ring element $\Theta_{S,T}(H/F, s) \in \C[G]$ satisfying \[ \chi(\Theta_{S,T}(H/F, s)) = L_{S,T}(\chi^{-1}, s) \quad \text{ for all } \chi \in \hat{G}.\]   
 Let $S_H$ denote the set of places of $H$ above those in $S$, and similarly for $T_H$.  Define
\[
U_{S,T} = \{x \in H^* : \text{ord}_w(x) \geq 0 \text{ for all } w \not\in S_H \text{ and } x \equiv 1 \!\!\pmod{T_H}\}.
\]
The condition on $T$ ensures that $U_{S,T}$ does not have any torsion. 
 The Galois equivariant version of Dirichlet's unit theorem 
 gives an $\R[G]$-module isomorphism 
\begin{equation} \label{e:lambda}
\lambda:  \R U_{S,T} \longrightarrow \R X_{S_H}
\end{equation}
\[
\lambda(u) = \sum_{w \in S_H} \log(|u|_w) \cdot w. 
\]

  Suppose that exactly $r$ places $v_1, \ldots, v_r \in S$ split completely in $H$ and $\# S \geq r+1$. The order of vanishing of $L_{S,T}(\chi, s)$ at $s=0$ is given by 
\[
r(\chi) = \dim_{\C} (\C U_{S,T})^{(\chi)} = \begin{cases}
 \#\{ v \in S\colon \chi(v) = 1\} \ \ \ \ \ & \text{ if } \chi \neq 1 \\ \#S - 1 & \text{ if } \chi = 1,
 \end{cases}
\]
whence $r(\chi) \ge r$ for all $\chi \in \hat{G}$.
Stark's conjecture predicts that the  $r$th derivative $\Theta_{S,T}^{(r)}(H/F,0)$ captures the ``non-rationality'' of the map $\lambda$.

\begin{conjecture}[Stark] \label{c:form1q} We have
\[
\Theta^{(r)}_{S,T}(H/F,0) \cdot \Q \bigwedge^r X_{S_H} \subset \lambda(\Q \bigwedge^r U_{S,T}).  
\]
\end{conjecture}

Concretely, this  states that for each character $\chi$ of $G$ with $r(\chi) = r$, the non-zero number $L^{(r)}_{S,T}(\chi^{-1}, s)$ lies in the one dimensional $\Q$-vector space spanned by $\lambda(\bigwedge^r (U_{S,T}^{(\chi)}))$.

Let us reformulate Conjecture~\ref{c:form1q} in terms of the existence of special elements. 
Write $X_{S_H}^* = \Hom(X_{S_H}, \mathbb{Z}[G])$. For $\varphi \in \bigwedge^r X_{S_H}^*$, there is a determinant map 
\[
\bigwedge^r X_{S_H} \times \bigwedge^r X_{S_H}^* \rightarrow \mathbb{Z}[G]
\] 
defined by 
\[
(x_1 \wedge \cdots \wedge x_r, \varphi_1 \wedge \cdots \wedge \varphi_r) \mapsto \varphi_1 \wedge \cdots \wedge \varphi_r(x_1 \wedge \cdots \wedge x_r) = \det(\varphi_i(x_j))_{i,j}.
\]  We extend the determinant map to $\R$-linearizations.
 We fix a place $w_i$ of $H$ above each $v_i$. Let $w_i^* \in X_{S_H}^*$ be induced by 
\[
w_i^*(w) = \sum_{\gamma \in G\colon \gamma w_i = w} \gamma.
\]

\begin{conjecture}[Stark] \label{c:form2q} Put $\varphi = w_1^* \wedge \cdots \wedge w_r^*$. There exists $u \in \Q \bigwedge^r U_{S,T}$ such that 
\[
\varphi(\lambda(u)) = \Theta^{(r)}_{S,T}(H/F,0).
\]

\end{conjecture}
The equivalence of Conjectures~\ref{c:form1q} and~\ref{c:form2q} is proven in \cite{rubin}*{Proposition 2.4}.

We are now ready to state the integral version of Stark's conjecture. 
In the rank $r=1$ case, Stark proposed the statement that $u$ in Conjecture~\ref{c:form2q} lies  in $U_{S,T}$.  This is the famous ``rank 1 abelian Stark conjecture.''
In the higher rank case,
the obvious generalization $u \in
\bigwedge^r U_{S,T}$ is not true, as was experimentally observed by Rubin \cite{rubin}. Rubin defined a lattice, nowadays called ``Rubin's lattice'' and  conjectured that it contains the element $u$. 

Put $U_{S,T}^* = \Hom_{\Z[G]}(U_{S,T}, \Z[G])$. 

The $r$th exterior bidual of $U_{S,T}$ (see \cite{burnssano} for a more general study and the initiation of this terminology) is defined by
\[
\bigcap^r U_{S,T} = 
\left(\bigwedge^r U_{S,T}^* \right)^* 
\cong \left\{ x \in \bigwedge^r \Q U_{S,T} : \varphi(x) \in \Z[G] \text{ for all } \varphi \in \bigwedge^r U_{S,T}^*\right\}.
\]

We would like to consider only the ``rank $r$'' component of this bidual.  To this end, for each character $\chi \in \hat{G}$ consider the  associated idempotent 
\[
e_{\chi} = \frac{1}{\#G} \sum_{g \in G} \chi(g) g^{-1} \in \C[G].
\]
Define
$e_r = \sum e_\chi \in \Q[G]$, where the sum extends over the set 
\[ \{\chi \in \hat{G}\colon L^{(r)}_{S,T}(\chi, 0) \neq 0\} =  \{\chi \in \hat{G}\colon \chi(G_v) \neq 1, v \in S \setminus \{v_1, \dots, v_r\}\}. \] 

Define Rubin's lattice by
\[
\mathcal{L}^{(r)} U_{S,T} =  \left( \bigcap^r U_{S,T}  \right) \cap e_r \left(\Q \bigwedge^r U_{S,T}^*\right).
\]

The following is Rubin's higher rank integral Stark conjecture.
\begin{conjecture}[\cite{rubin}, Conjecture B'] \label{c:form1z} Put $\varphi = w_1^* \wedge \cdots \wedge w_r^*$. There exists $u \in  \mathcal{L}^{(r)} U_{S,T}$ such that 
\[
\varphi(\lambda(u)) = \Theta^{(r)}_{S,T}(H/F,0).
\]
\end{conjecture}

\section{Stark's Conjectures at finite places} \label{s:starkfinite}

 We now assume that the totally split places $v_1, \ldots, v_r$ from the previous section are all finite. 
  This happens only when $F$ is a totally real field and $H$ is totally complex.
  In fact, the fixed fields of characters with nonvanishing $L$-functions at 0 are CM fields, so we restrict to the setting where $F$ is totally real and $H$ is CM for the remainder of the article.
 We also enact a slight notational change and write the set denoted $S$ in the previous sections as $S'$, and let $S = S' \setminus \{v_1, \dots, v_r\}$. The reason for this is that  we now still have $S \supset S_\infty \cup S_{\ram}$. 

 As we explain, in this setting Conjecture \ref{c:form2q} for $S'$ follows from a classical rationality result of Klingen--Siegel, though the integral refinement in Conjecture~\ref{c:form1z} remains a nontrivial statement. For a fixed place $w$ of $H$, we have
\begin{equation} \label{e:lo}
\log |u|_w = - \ord_w(u) \log \N w.
\end{equation}
Since the Euler factors at the $v_i$ are equal to $(1 - \N v_i^{-s}) = (1 - \N w_i^{-s})$, we also have
\begin{equation} \label{e:thetarel}
\Theta_{S',T}^{(r)}(H/F,0) = \Theta_{S,T}(H/F,0) \cdot \prod_{i=1}^r \log \N w_i.
\end{equation}
\begin{theorem}[Klingen--Siegel]  \label{t:ks}
 We have $\Theta_{S,T} := \Theta_{S,T}(H/F,0) \in \Q[G]$.
\end{theorem}

With $e_r$ as in the previous section, we are then led to define a map over $\Q$
\[
\lambda_{\Q} : e_r(\Q U_{S',T}) \rightarrow e_r(\Q X_{S'_H}), \qquad
\lambda_{\Q}(u) = \sum_{w | v_i \text{ some } i} \ord_w(u) \cdot w.
\]
Note that $ e_r(\Q X_{S'_H})$ is the $\Q$-vector space generated by the places of $H$ above the $v_i$.
The map $\lambda_{\Q}$ is a $\Q[G]$-module isomorphism, and it induces an isomorphism on the free rank one $\Q[G]$-modules obtained by taking $r$th wedge powers.  In view of (\ref{e:lo}), the map on $r$th wedge powers induced by the map $\lambda$ of (\ref{e:lambda}), when
restricted to $e_r(\Q \bigwedge^r U_{S',T})$, is equal to $(\prod_{i=1}^r \log \N w_i) \cdot \lambda_{\Q}$.
Conjecture~\ref{c:form2q}  follows  from this observation together with (\ref{e:thetarel}), since  Theorem~\ref{t:ks}
implies the existence of $u \in e_r(\Q \bigwedge^r U_{S',T})$ such that 
\[
\varphi(\lambda_{\Q}(u)) = \Theta_{S,T}.
\]
Here $\varphi = w_1^* \wedge \cdots \wedge w_r^*$ as in the statement of the conjecture.

On the other hand, the integral statement in Conjecture \ref{c:form1z} lies deeper.  We first note the following celebrated theorem of Deligne--Ribet and Cassou-Nogu\`es refining the Klingen--Siegel theorem. The condition on the set $T$ stated in \S \ref{s:stark} is crucial in this result. 
\begin{theorem} We have $\Theta_{S,T} \in \Z[G]$.
\end{theorem}

Conjecture \ref{c:form1z} in this setting is known as the Rubin--Brumer--Stark conjecture:
\begin{conjecture}[Rubin--Brumer--Stark] \label{c:form2z} There exists $u \in  \mathcal{L}^{(r)} U_{S',T}$ such that 
\[
\varphi(\lambda_{\Q}(u)) = \Theta_{S,T}.
\]
\end{conjecture}
We describe in Theorem~\ref{t:ikur} below a strong partial result toward this conjecture.

\section{The Brumer--Stark Conjecture}  \label{s:bs}

Having stated the higher rank Rubin--Brumer--Stark conjecture, we now wind back the clock and focus on the case $r=1$.  This case had been studied independently by Brumer and Stark and served as a motivation for Rubin's work.  Writing the splitting prime $v_1$ as $\fp$, the conjecture may be stated as follows.
\begin{conjecture}[Brumer--Stark] \label{c:bs} Fix a prime ideal $\fp \subset O_F, \fp \not\in S \cup T$, such that $\fp$ splits completely in $H$.
 Fix a prime $\fP \subset O_H$ above $\fp$. There exists a unique element $u_\fp \in H^*$ such that $|u_\fp|_w=1$ for every place $w$ of $H$ not lying above $\fp$, 
\begin{equation} \label{e:bs}
 \ord_G(u_\fp) := \sum_{\sigma \in G} \ord_{\fP}(\sigma(u_\fp)) \sigma^{-1} = \Theta_{S,T}, \end{equation}
and $u \equiv 1 \pmod{\fq}$ for all $\fq \in T_H$.
\end{conjecture}
Note that the condition $|u|_w =1$  includes all complex places $w$, so $c(u_\fp) = u_\fp^{-1}$ for the unique complex conjugation $c \in G$.

As we have alread noted, Stark arrived upon this statement in the 1970s through his attempts to generalize and factorize the classical Dirichlet class number formula (though in a slightly different formulation; the statement above is due to Tate \cite{tate}). Prior to this, in the 1960s,
Brumer was interested in generalizing Stickelberger's classical factorization formula for Gauss sums in cyclotomic fields.  Stickelberger's result can be formulated as stating that when $H = \Q(\mu_N)$ is a cyclotomic field, the Stickelberger element annihilates the class group of $H$.  Let us consider Brumer's perspective of annihilation of class groups  in the case of general $H/F$.

\subsection{Annihilation of Class groups}

Let $\Cl^T(H)$ denote the ray class group of $H$ with conductor equal to the product of primes  in $T_H$.  This is defined as follows.
Let $I_T(H)$ denote the group of fractional ideals of $H$ relatively prime to the primes in $T_H$.
  Let $P_{T}(H)$ denote the subgroup of $I_T(H)$ generated by principal ideals $(\alpha)$ where $\alpha \in O_H$ satisfies $\alpha \equiv 1 \pmod{\fq}$ for all $\fq \in T_H$.  Then \[ \Cl^T(H) = I_T(H)/P_{T}(H). \]  This $T$-smoothed class group is naturally a $G$-module. 
  
  With the notation as in Conjecture~\ref{c:bs}, we have \begin{equation} \label{e:fp} 
  \fP^{\Theta_{S,T}} = (u_\fp). \end{equation} Such an equation holds for all $\fp \not\in S \cup T$ that split completely in $H$. The set of primes of $H$ above all such $\fp$ generate $\Cl^T(H)$. Hence we deduce
\begin{equation} \label{e:bs2}
  \Theta_{S,T} \in \Ann_{\Z[G]}(\Cl^T(H)).  \end{equation}
In fact (\ref{e:bs2}) is almost equivalent to Conjecture \ref{c:bs}; given (\ref{e:fp}), the element $u_\fp$ satisfies all of the conditions necessary for Conjecture \ref{c:bs} except possibly $c(u_\fp) = u_\fp^{-1}$.  But of course $v_{\fp} = u_\fp/c(u_\fp)$ satisfies this condition and moreover satisfies $\fP^{2\Theta_{S,T}} = (v_\fp)$.  Therefore the only possible discrepancy between the statements is a factor of 2, which disappears when we localize away from 2 as in the rest of this paper.  Let us therefore define \[ R = \Z[1/2][G]^- = \Z[1/2][G]/(c+1), \]
and for any $\Z[G]$-module $M$ we write $M^- = M \otimes_{\Z[G]} R$.  There exists an element $u_\fp \in O_H[1/\fp]^* \otimes \Z[1/2]$ satsifying Conjecture~\ref{c:bs} if and only if 
\begin{equation} \label{e:bs3}
  \Theta_{S,T} \in \Ann_{R}(\Cl^T(H)^-).  \end{equation}
  This is the Brumer--Stark conjecture ``away from 2.''
  
Many authors have studied (\ref{e:bs3}) as well as refinements.  The works of Burns, Greither, Kurihara, Popescu, and Sano are particularly noteworthy (\cites{greitherkurihara, greither, kurihara, burns1, burns, burnssano, greitherpopescuemc}).  For instance, it is very natural to ask whether 
$ \Theta_{S,T}$  lies in the 0th Fitting ideal of $\Cl^T(H)^-$ over $R$, since the Fitting ideal is contained in the annihilator.  (See the beginning of \S\ref{s:refine} to recall the definition of Fitting ideal.)  It was  noticed by Popescu in the function field case that while this holds sometimes, it does not always hold \cite{popfunc}.  
Upon hearing of these examples, Kurihara found similar examples of non-containment in the number field setting.
 Greither and Kurihara (\cites{greither,greitherkurihara}) observed that the statement may be corrected by replacing $\Cl^T(H)^-$ by its Pontryagin dual \[ \Cl^T(H)^{-,\vee} = \Hom_{\Z}(\Cl^T(H)^-, \Q/\Z). \]
We endow $\Cl^T(H)^{-,\vee}$ with the contragradient $G$-action $g\cdot \varphi(x) = \varphi(g^{-1}x)$.  Denote by $\#$ the involution on $\Z[G]$ induced by $g \mapsto g^{-1}$ for $g \in G$.
\begin{conjecture}[Kurihara, ``Strong Brumer--Stark"] \label{c:sbs}
We have
\[ \Theta_{S,T}^\# \in \Fitt_R(\Cl^T(H)^{-,\vee}). \]
\end{conjecture}
Conjecture~\ref{c:sbs} leads to the following natural questions.
\begin{enumerate}
\item What is the Fitting ideal of $\Cl^T(H)^{-, \vee}$?
\item What is the Fitting ideal of $\Cl^T(H)^-$?
\item Is there a natural arithmetically defined $R$-module whose Fitting ideal is generated by $\Theta_{S,T}$ or $\Theta_{S,T}^\#$?
\end{enumerate}

The precise conjectural description of the Fitting ideal of $\Cl^T(H/F)^{-,\vee}$ was given by Kurihara \cite{kurihara}; we state this in \S\ref{s:kur} below.  An important fact about this statement is that when $S_{\ram}$ is nonempty, the 
Fitting ideal of $\Cl^T(H/F)^{-,\vee}$ is in general not principal (and in particular is not generated by $\Theta_{S,T}^\#$).

A conjectural answer to the second question above has recently been provided in a striking paper by Atsuta and Kataoka \cite{atsutakataoka}. They show that their conjecture is implied by the Equivariant Tamagawa Number Conjecture. 

The third question is answered by a conjecture of Burns and is the topic of \S\ref{s:burns}.  We note that Fitting ideals of finitely presented $R$-modules are rarely principal. It is therefore remarkable that Burns defined a natural arithmetic object whose Fitting ideal is principal. 

\subsection{Our Results}
We now describe some of our results toward these conjectures \cite{bsdk}*{Theorem 1.4}.
\begin{theorem}  \label{t:ikur}
Kurihara's exact formula for $\Fitt_R(\Cl^T(H)^{-,\vee})$ holds (see Theorem~\ref{t:kur}).  In particular, we have the Brumer--Stark and Strong Brumer--Stark conjectures away from 2:
 \[ \Theta_{S,T}^\# \in \Fitt_R(\Cl^T(H)^{-,\vee}) \subset \Ann_R(\Cl^T(H)^-)^\#. 
\]
Finally, Rubin's higher rank Brumer--Stark conjecture holds away from 2: with notation as in Conjecture~\ref{c:form2z}, 
there exists $u \in  \mathcal{L}^{(r)} U_{S',T} \otimes \Z[1/2]$ such that 
$\varphi(\lambda_{\Q}(u)) = \Theta_{S,T}.$

\end{theorem}

Partial results in this direction had been obtained earlier by Burns  \cite{burns}  (including a $\mu=0$ hypothesis and the assumption of the Gross--Kuz'min conjecture) and by  Greither and Popescu  \cite{greitherpopescuemc} (including a $\mu=0$ hypothesis and imprimitivity conditions on $S$).

Our results in \cite{bsdk} do not seem to directly imply the conjecture of Atsuta and Kataoka on $\Fitt_R(\Cl^T(H)^-)$ or the conjecture of Burns. 
However we prove an analogous result toward the latter, with $(S,T)$ replaced by an alternate pair $(\Sigma, \Sigma')$, in Theorem~\ref{t:selfit}.  This result turns out to be strong enough to deduce Theorem~\ref{t:ikur}.

In \S\ref{s:rw} we give a detailed summary of the proof of Theorem~\ref{t:selfit}.  Key ingredients are the $\Z[G]$-modules $\nabla_{\Sigma}^{\Sigma'}(H)$ defined by Ritter and Weiss, and Ribet's method of  using modular forms to construct Galois cohomology classes associated to $L$-functions.

\subsection{Explicit Formula for Brumer--Stark units} \label{s:ef}

We conclude this section by describing a further direction in the study of Brumer--Stark units, that of {\em explicit formulae} and applications to explicit class field theory.  This theme was
initiated by Gross in \cite{gross} and \cite{gross2} and developed by the first author and collaborators over a series of papers (\cites{darmondasgupta, dasguptaduke, charolloisdasgupta, cdg, dasguptaspiess}).

Let $\fp$ be as above and write $S' = S \cup \{\fp\}$.
Let $L$ denote a finite abelian CM extension of $F$ containing $H$ that is ramified over $F$ only at the places in $S'$.  Write $\fg = \Gal(L/F)$ and $\Gamma = \Gal(L/H)$, so $\fg/\Gamma \cong G$.  Let $I$ denote the relative augmentation ideal associated to $\fg$ and $G$, i.e.\ the kernel of the canonical projection $
 \Aug\colon \Z[ \fg] \longtwoheadrightarrow \Z[G].$  Then $ \Theta_{S',T}({L/F})   $ lies in $I$, since its image under $\Aug$ is 
\begin{equation} \label{e:euler}
 \Theta_{S',T}({H/F})= \Theta_{S,T}({H/F}) (1 - \sigma_\fp) = 0. \end{equation}
Here $\sigma_\fp$ denotes the Frobenius associated to $\fp$ in $G$, and this is trivial since
 $\fp$ splits completely in $H$.  Intuitively, if we view $ \Theta_{S',T}(L/F)$ as a function on the ideals of $\Z[\fg]$, equation (\ref{e:euler}) states that this function ``has a zero'' at the ideal $I$; the value of the ``derivative'' of this function at $I$ is simply the image of $ \Theta_{S',T}(L/F) $ in $I/I^2$.  Gross provided a conjectural algebraic interpretation of this derivative as follows.
Denote by \[
 \rec_\fP \colon H_\fP^* \longrightarrow \Gamma\]
 the composition of the inclusion $H_\fP^* \longhookrightarrow \A_H^*$ with the global Artin reciprocity map \[ \A_H^* \longtwoheadrightarrow \Gamma. \]
Throughout this article we adopt Serre's convention \cite{serre} for the reciprocity map. Therefore $\rec(\varpi^{-1})$ is a lifting to $G_{\fp}^{\ab}$ of the Frobenius element on the maximal unramified extension of $F_{\fp}$ if $\varpi \in F_{\fp}^*$ is a uniformizer. 

\begin{conjecture}[Gross, \cite{gross2}*{Conjecture 7.6}] \label{c:grosstower} Define
\begin{equation} \label{e:recg}
\rec_G(u_\fp) = \sum_{\sigma \in G} (\rec_\fP\sigma(u_\fp)-1)  \tilde{\sigma}^{-1} \in I/I^2, \end{equation}
where $\tilde{\sigma} \in \fg$ is any lift of $\sigma \in G$.  Then
\[ \rec_G(u_\fp) \equiv  \Theta_{S',T}^{L/F}  \quad \text{ in } I/I^2. \]
\end{conjecture}

The main result of \cite{hdk} is  the following.

\begin{theorem} \label{t:maingross} Let $p$ be an odd prime and suppose that $\fp$ lies above $p$. 
Gross's Conjecture~\ref{c:grosstower} holds in $(I/I^2) \otimes \Z_p$.
 \end{theorem}
 
 Our interest in this result is that by enlarging $S$ and taking larger and larger field extensions $L/F$, one can use (\ref{e:recg}) to specify all of the $\fp$-adic digits of $u_{\fp}$.
One therefore obtains an exact $\fp$-adic analytic formula for $u_\fp$.  This formula can be described  either using the {\em Eisenstein cocycle} or more explicitly via {\em Shintani's method}; the latter approach is followed in \S\ref{s:exact}.
 In \S\ref{s:hit}, we describe the argument using ``horizontal Iwasawa theory'' to show that Theorem~\ref{t:maingross} implies the conjectural exact formula.  In \S\ref{s:gs} we summarize the key ingredients involved in the proof of Theorem~\ref{t:maingross}, including an integral version of the Greenberg--Stevens $\sL$-invariant and an associated modified Ritter--Weiss module $\nabla_\sL$.
In the setting of $F$ real quadratic, Darmon, Pozzi, and Vonk have given an alternate, elegant proof of the explicit formula for the units $u_\fp$ (\S\ref{s:dpv}).  Their approach involves $p$-adic deformations of modular forms, rather than the tame deformations that we consider.

One significance of the exact formula is that we show that the collection of Brumer--Stark units, together with some easily described elements, generate the maximal abelian extension of the totally real field $F$.
  \begin{theorem}  \label{t:h12}
  Let $\BS$ denote the set of Brumer--Stark units $u_\fp$ as we range over all possible CM abelian extensions $H/F$ and for each extension a choice of prime $\fp$ that splits completely in $H$.   Let $\{\alpha_1, \dotsc, \alpha_{n-1}\}$ denote any elements of $F^*$ whose signs in $\{\pm 1\}^{n}/ (-1, \dotsc, -1)$ under the real embeddings of $F$ form a basis for this $\Z/2\Z$-vector space.  
The maximal abelian extension of $F$ is generated by $\BS$ together with $\sqrt{\alpha_1}, \dotsc, \sqrt{\alpha_{n-1}}$:
\[ F^{\ab}  =  \ F(\BS, \sqrt{\alpha_1}, \dotsc, \sqrt{\alpha_{n-1}}). \]
\end{theorem}
It is important to stress that the exact formula for $u_\fp$ described in \S\ref{s:exact}
can be computed without knowledge of the field $H$,  using only the data of $F$, $\fp$, and the conductor of $H/F$. Furthermore, we can leave out any finite set of primes $\mathfrak{p}$ without altering the conclusion of the theorem. In this way we obtain an {\em explicit class field theory} for $F$, i.e.\ an analytic construction of its maximal abelian extension $F^{\ab}$ using data intrinsic only to $F$ itself. Explicit computations of class fields of real quadratic fields generated using our formula are provided in \cite{bsdk}*{\S2.3} and \cite{fl}.

\section{Refinements of Stark's conjecture} \label{s:refine} In this section we recall various refinements of the strong Brumer--Stark conjecture.   We first recall the definition of Fitting ideal.
 Let $R$ be a commutative ring and $M$  an $R$-module with finite presentation:
\[
R^m \xrightarrow{A} R^n \rightarrow X \rightarrow 0.
\]
Here $A$ is an $n \times m$ matrix over $R$.  
The $i$th Fitting ideal $\Fitt_{i, R}(M)$ is the ideal generated by the $n-i \times n-i$ minors of $A$. It is a standard fact \cite{northcott}*{Chapter 3, Theorem 1} that 
$\Fitt_{i,R}(M)$ does not depend on the chosen presentation of $M$.  We  write $\Fitt_{R}(M)$ for $\Fitt_{0,R}(M)$.

\subsection{The conjecture of Kurihara} \label{s:kur}
 In this section we describe the Fitting ideal of the minus part of the dual class group. 
For each $v$ in $S_{\ram}$, let $I_v \subset G_v \subset G$ denote the inertia and decomposition groups, respectively, associated to $v$. Write
\[
e_v = \frac{1}{\# I_v} \N I_v = \frac{1}{\# I_v} \sum_{\sigma \in I_v} \sigma \in \Q[G]
\]
for the idempotent that represents projection onto the characters of $G$ unramified at $v$. Let $\sigma_v \in G_v$ denote any representative of the Frobenius coset of $v$. The element $1 - \sigma_v e_v \in \Q[G]$ is independent of choice of representative. Following \cite{greither}, we define the Sinnott--Kurihara ideal, \emph{a priori} a fractional ideal of $\Z[G]$, by 
\[
{\Sku}^T(H/F) = (\Theta^{\#}_{S_{\infty}, T}) \prod_{v \in S_{\ram}}(\N I_v, 1 - \sigma_v e_v).
\]
Kurihara proved using the theorem of Deligne--Ribet and Cassou-Nogu\`es that $\Sku^T(H/F)$ is a subset of $\Z[G]$ (see \cite{bsdk}*{Lemma 3.4}). The following conjecture of Kurihara is proven in \cite{bsdk}*{Theorem 1.4}.

\begin{theorem} \label{t:kur} We have
\[
\Fitt_{R}(\Cl^T(H)^{-,\vee}) = \Sku^T(H/F)^-.
\]
\end{theorem}

 The definition of the Sinnott--Kurihara ideal is inspired by Sinnott's definition of generalized Stickelberger elements for abelian extensions of $\Q$  \cite{sinnott}. For a generalization of Sinnott's ideal to arbitrary totally real fields see \cite{greither}*{\S2}. In general, Sinnott's ideal contains the Sinnott--Kurihara ideal but it may be strictly larger.  

The plus part of the Sinnott--Kurihara ideal is not very interesting as the plus part of $\Theta^{\#}_{S_\infty, T}$ is 0. The plus part of the class group is much smaller than the minus part and seems  harder to describe; for example Greenberg's conjecture on the vanishing of lambda invariants implies that the order of the plus part is bounded up the cyclotomic tower. For abelian extensions of $\Q$, the plus part is described by Sinnott using  cyclotomic units.

\subsection{The conjecture of Atsuta--Kataoka} It is in fact more natural to ask about the Fitting ideal of $\Cl^T(H)$, as opposed to the Pontryagin dual. A conjectural answer to this question has been provided in a recent paper of Atsuta--Kataoka \cite{atsutakataoka} using the theory of shifted Fitting ideals developed by Kataoka \cite{kataoka}. 
We recall this notion now. 

Let $M$ be an $R$-module of finite length. Take a resolution 
\[
0 \rightarrow N \rightarrow P_1 \rightarrow \cdots \rightarrow P_d \rightarrow M \rightarrow 0
\]
with each $P_i$ of projective dimension $\leq$ 1. Following \cite{kataoka} define the shifted Fitting ideal
\[
\Fitt_R^{[d]}(M) = \left( \prod_{i=1}^d \Fitt_R(P_i)^{(-1)^i} \right) \Fitt_R(N).
\]
The independence of this definition from the choice of resolution is proven in \cite{kataoka}*{Theorem 2.6 and Proposition 2.7}. 
Let \[
g_v = 1 - \sigma_v + \#I_v \in \Z[G/I_v], \qquad h_v = 1 - e_v \sigma_v + \N I_v \in \Q[G].
\]
Define the $\Z[G]$-module 
\[
A_v = \Z[G/I_v]/(g_v).
\]

\begin{conjecture}[Atsuta--Kataoka] \label{c:ak} We have
\[
\Fitt_{R}(\Cl^T(H)^-) = \left( \prod_{w \in S_{\ram, H}} h_v^- \Fitt_{R}^{[1]}(A_v^{-}) \right) \Theta_{S_{\infty}, T}.
\]
\end{conjecture}

In \cite{atsutakataoka}, the authors give an explicit description of the ideal $h_v^- \Fitt_{\Z[G]^-}^{[1]}(A_v^-)$ appearing in Conjecture~\ref{c:ak}. 
Write $I_v = J_1 \times \cdots \times J_s$ for cyclic groups $J_i$, $1 \leq i \leq s$. For each $i$  put 
\[
\N_i = \N J_i = \sum_{\sigma \in J_i} \sigma \in \Z[G].
\] 
Furthermore, put $\mathcal{I} = \ker(\Z[G] \rightarrow \Z[G/G_v])$ for the relative augmentation ideal. For each $1 \leq i \leq s$, put $Z_i$ for the ideal of $\Z[G]$ generated by $\N_{j_1} \cdots \N_{j_{s-i}}$, where $(j_1, \ldots, j_{s-i})$ runs through all tuples of integers satisfying $1 \leq j_1 \leq \cdots \leq j_{s-i} \leq s$. This definition of $Z_i$ is independent of the decomposition of $I_v$ into cyclic groups. Define
\[
\mathcal{J} = \sum_{i=1}^s Z_i \mathcal{I}^{i-1}.
\]

\begin{theorem}[Atsuta--Kataoka] We have
\[
h_v^- \Fitt_{\Z[G]^-}^{[1]}(A_v^-) = \left( \N I_v, \left( 1 - e_v \sigma_v \right) \mathcal{J} \right)
\]
as fractional ideals of $\Z[G]^-$.
\end{theorem}
Atsuta--Kataoka prove:
\begin{theorem} The Equivariant Tamagawa Number Conjecture for $H/F$ implies Conjecture~\ref{c:ak}.
\end{theorem}

\subsection{The conjecture of Burns} \label{s:burns}
The refinements mentioned above do not involve principal ideals. The method of Ribet, which attempts to show the inclusion of an arithmetically defined ideal into an analytically defined ideal, works well for principal ideals. From this point of view it is   natural to ask if there is an arithmetically defined object whose Fitting ideal is generated by the Stickelberger element 
$\Theta_{S,T}$. Burns provided a conjectural answer to this question \cite{burns}.  A modification of this statement (Theorem \ref{t:selfit} below) is the main technical result in \cite{bsdk} from which all the results of Theorem~\ref{t:ikur} are deduced. 

We now recall the statement of Burns's conjecture.
Let $H^*_{T}$ be the group of $x \in H^*$ such that $\ord_{w}(x-1) > 0$ for each prime $w \in T_H$. Define
\[
\Sel^{T}_S(H) = \Hom_{\Z}(H_{T}^*, \Z)/ \prod_{ w \notin S_H \cup {T}_H} \Z,
\]
where the implicit map sends a tuple $(x_w)$ to the function $\sum_w x_w \ord_w$. The $G$-action on $\Sel^{T}_S(H)$ is the contragradient $G$-action $(g \varphi)(x) = \varphi(g^{-1}x)$. 

\begin{conjecture}[Burns] We have 
\[
\Fitt_{R}(\Sel^T_S(H)^-) = (\Theta_{S,T}^\#).
\]
\end{conjecture}

We have proven a version of this result with altered sets $S$ and $T$.  Fix an odd prime $p$ and put $R_p = \Z_p[G]^-$.  Define
\[
\Sigma = S \setminus \{ v \in S : v \nmid p\}
\]
and 
\[
\Sigma' = T \cup \{v \in S : v \nmid p\}
\]
\begin{theorem}[\cite{bsdk}, Theorem 3.3] \label{t:selfit}
Let $\Sel^{\Sigma'}_{\Sigma}(H)_p^- = \Sel^{\Sigma'}_{\Sigma}(H) \otimes_{\Z[G]} R_p$.  We have
\[
\Fitt_{R_p}(\Sel^{\Sigma'}_{\Sigma}(H)_p^-) = (\Theta_{\Sigma, \Sigma'}^{\#}).
\]
\end{theorem}
It turns out that Theorem~\ref{t:selfit} is strong enough to imply Kurihara's conjecture (Theorem~\ref{t:kur}).
The key point is that there is a short exact sequence 
\begin{equation} \label{e:sst0}
\begin{tikzcd}
 0 \ar[r] &   \Sel_{\Sigma}^T(H)^-  \ar[r] & \Sel_{\Sigma}^{\Sigma'}(H)^- \ar[r] & \prod_{w \in S'_H} ((O_{H}/w)^*)^{\vee, -}  \ar[r] & 0,
 \end{tikzcd}
\end{equation}
from which one deduces (see \cite{bsdk}*{Theorem 3.7}) 
\[  \Fitt_{R_p}(\Sel_{\Sigma}^T(H)_p^-) = (\Theta^\#_{\Sigma,T}) \prod_{v \in S_{\ram},\ \! v\nmid p} \!\!\!\! ( \N I_v, 1 - \sigma_v e_v ). \]
Since $\Sel_{S_\infty}^T(H)^- \cong \Cl^T(H)^{-, \vee}$, it then remains to calculate the effect of removing the primes $v \in S_{\ram}, v \mid p$ from $\Sigma$.  This is a delicate process using functorial properties of the Ritter--Weiss modules discussed in \S\ref{s:rw}, and one obtains (see \cite{bsdk}*{Appendix B}) the desired result
\[  \Fitt_{R_p}(\Sel_{S_\infty}^T(H)_p^-) = (\Theta^\#_{S_\infty,T}) \prod_{v \in S_{\ram}}  ( \N I_v, 1 - \sigma_v e_v ). \]

\section{Ritter--Weiss modules and Ribet's method}  \label{s:rw}

In this section we summarize the proof of Theorem~\ref{t:selfit}.

\subsection{Ritter--Weiss modules}

The $\Z[G]$-module that shows up in our constructions with modular forms is a certain {\em transpose} of $\Sel_S^T(H)$ in the sense of Jannsen \cite{jannsen}, denoted $\nabla_S^T(H)$.
This module was originally defined by Ritter and Weiss in the foundational paper~\cite{rw} without the smoothing set $T$.  We incorporated the smoothing set $T$ and established some additional properties of $\nabla_S^T(H)$ in \cite{bsdk}*{Appendix A}.  To describe these properties, we work over $R_p=\Z_p[G]^-$ and consider finite disjoint sets $\Sigma, \Sigma'$ satisfying the following:

\begin{itemize} \item $\Sigma \supset S_\infty$ and  $\Sigma \cup \Sigma' \supset S_{\ram}$
 \item $\Sigma'$ satisfies the condition (\ref{i:conditions}) on $T$ in \S \ref{s:stark}.
 \item The primes in $\Sigma' \cap S_{\ram}$ have residue characteristic $\ell \neq p$.
\end{itemize}
Note that the pair $(S,T)$ from \S \ref{s:stark} and the pair $(\Sigma, \Sigma')$ considered in \S\ref{s:burns} both satisfy these conditions.
The module $\nabla_{\Sigma}^{\Sigma'}(H)_p^- = \nabla_{\Sigma}^{\Sigma'}(H) \otimes_{\Z[G]} R_p$ satisfies the following:
\begin{itemize}
\item There is a short exact sequence of $R_p$-modules
\begin{equation} \label{e:nablaext}
 \begin{tikzcd}
 0 \ar[r] &  \Cl_\Sigma^{\Sigma'}(H)_p^- \ar[r] &  \nabla_{\Sigma}^{\Sigma'}(H)_p^- \ar[r] &  (X_{H_\Sigma})_p^- \ar[r] &  0.
 \end{tikzcd}
  \end{equation}
Here $\Cl_\Sigma^{\Sigma'}(H)$ denotes the quotient of $\Cl^{\Sigma'}(H)$ by the image of the primes in $\Sigma_H$.

 \item Given a $R_p$-module $B$, a surjective $R_p$-module homomorphism \begin{equation} \label{e:nablasurj}
  \nabla_{\Sigma}^{\Sigma'}(H)_p^- \longtwoheadrightarrow B
  \end{equation} is equivalent to the data of a cocycle $\kappa \in Z^1(G_F, B)$ and a collection of elements $x_v \in B$ for $v \in \Sigma$ satsifying the following conditions:
 \begin{itemize}
  \item The cohomology class $[\kappa] \in H^1(G_F, B)$ is unramified outside $\Sigma'$, tamely ramified at $\Sigma'$, and locally trivial at $\Sigma$.
  \item The $x_v$ provide  local trivializations at $\Sigma$: $\kappa(\sigma) = (\sigma - 1)x_v$ for $\sigma \in G_v$. 
\item The $x_v$ along with the image of $\kappa$ generate the module $B$ over $R_p$.
\end{itemize}
   The tuples $(\kappa, \{x_v\})$ are taken modulo the natural notion of coboundary, i.e.\ $(\kappa, \{x_v\}) \sim (\kappa + dx, \{x_v + x\})$ for $x \in B$.
\item The module $\nabla_{\Sigma}^{\Sigma'}(H)_p^-$ has a quadratic presentation, i.e.\ there exists an exact sequence of $R_p$-modules
\begin{equation} \label{e:qp}
 \begin{tikzcd} M_1 \ar[r,"A"] & M_2 \ar[r] &   \nabla_{\Sigma}^{\Sigma'}(H)_p^-  \ar[r] &  0 \end{tikzcd}
\end{equation}
where $M_1$ and $M_2$ are free of the same finite rank.
\item The module $ \nabla_{\Sigma}^{\Sigma'}(H)_p^-$ is a transpose of $\Sel_{\Sigma}^{\Sigma'}(H)_p^-$, i.e. for a suitable quadratic presentation (\ref{e:qp}) of $ \nabla_{\Sigma}^{\Sigma'}(H)_p^-$, the cokernel of the induced map
\begin{equation} \label{e:tr}
  \begin{tikzcd} \Hom_{R_p}(M_2, R_p) \ar[r,"A^{T, \#}"] & \Hom_{R_p}(M_1, R_p) \end{tikzcd}
\end{equation}
is isomorphic to $\Sel_{\Sigma}^{\Sigma'}(H)_p^-$.  Here we follow our convention of giving  Hom spaces the contragradient $G$-action.
\end{itemize}

The quadratic presentation property (\ref{e:qp}) implies that $\Fitt_{R_p}(\nabla_{\Sigma}^{\Sigma'}(H)_p^-) = \det(A)$ is principal.  The transpose property (\ref{e:tr}) implies that
\begin{equation} \label{e:pound}
 \Fitt_{R_p}(\nabla_{\Sigma}^{\Sigma'}(H)_p^-) = \Fitt_{R_p}(\Sel_{\Sigma}^{\Sigma'}(H)_p^-)^{\#}. \end{equation}
Theorem~\ref{t:selfit} is therefore equivalent to  \begin{equation} \label{e:fittnabla}
 \Fitt_{R_p}(\nabla_{\Sigma}^{\Sigma'}(H)_p^-) = (\Theta_{\Sigma, \Sigma'}).
 \end{equation}
 We now fix $(\Sigma, \Sigma')$ to be the pair defined in \S\ref{s:burns}.
In the remainder of this section we outline how (\ref{e:fittnabla}) is proved using Ribet's method.
Throughout, an unadorned $\Theta$ denotes $\Theta_{\Sigma, \Sigma'}$ (and $\Theta^\#$ denotes $\Theta_{\Sigma, \Sigma'}^\#$).

\subsection{Inclusion Implies Equality}

An interesting feature of Ribet's method is that it tends to prove an inclusion in one direction, that of an algebraically defined ideal contained within an analytically defined ideal. In our setting, we use it to prove
\begin{equation} \label{e:inclusion}
  \Fitt_{R_p}(\nabla_{\Sigma}^{\Sigma'}(H)_p^-) \subset (\Theta), \quad \text{equivalently,} \quad
  \Fitt_{R_p}(\Sel_{\Sigma}^{\Sigma'}(H)_p^-) \subset (\Theta^\#).
   \end{equation}
We then employ an analytic argument to show that this inclusion is an equality.  It is important to note that the inclusion (\ref{e:inclusion}) is the {\em reverse direction} of that required by the Brumer--Stark and Strong Brumer--Stark conjectures.  It is therefore essential for our approach that one actually has the statement of an {\em equality} rather than just an inclusion (and an analytic argument to deduce the equality from the reverse inclusion).
For this reason, the conjecture of Burns stated in \S\ref{s:burns} (more precisely the analog of it stated in Theorem~\ref{t:selfit}) plays an essential role in our strategy.

Let us describe the analytic argument in the special case  $\Sigma = S_\infty$, i.e.\ there are no primes above $p$ ramified in $H/F$.  
In this case 
\begin{equation} \label{e:isomsc}
\Sel_{\Sigma}^{\Sigma'}(H)_p^- \cong \Cl^{\Sigma'}(H)_p^- \end{equation} is finite and $\Theta$ is a non-zerodivisor.
Using (\ref{e:inclusion}), write \begin{equation} \label{e:selz}
 \Fitt_{R_p}(\Sel_{\Sigma}^{\Sigma'}(H)_p^-) = (\Theta^{\#} \cdot z) \quad \text{ for some } z \in R_p. \end{equation} We must show that $z \in R_p^*$.
An elementary argument (see \cite{bsdk}*{\S 2.3})  shows that (\ref{e:selz}) implies
\begin{equation} \label{e:selsize}
 \# \Sel_{\Sigma}^{\Sigma'}(H)_p^- = \bigg(\prod_{\substack{\chi \in \hat{G} \\ \chi \text{ odd}}}  \chi(\Theta^{\#} \cdot z)  \bigg)_p
\end{equation}
where the subscript $p$ on the right denotes the $p$-power part of an integer.
Yet the analytic class number formula implies (see \cite{bsdk}*{\S2.1}) 
\begin{equation} \label{e:prod}  \prod_{\substack{\chi \in \hat{G} \\ \chi \text{ odd}}}  \chi(\Theta^{\#}) =  
\prod_{\substack{\chi \in \hat{G} \\ \chi \text{ odd}}} L_{\Sigma, \Sigma'}(\chi, 0) \doteq \# \Cl^{\Sigma'}(H)^-, 
\end{equation}
where $\doteq$ denotes equality up to a power of 2.  Combining (\ref{e:isomsc}), (\ref{e:selsize}), and (\ref{e:prod}), one finds that $\chi(z)$ is a $p$-adic unit for each odd character $\chi$.  It follows that $z \in R_p^*$ as desired.

The generalization of this argument to arbitrary $\Sigma$ requires a delicate induction and is described in \cite{bsdk}*{\S5}.

\subsection{Ribet's method}

We now describe our implementation of Ribet's method to prove the inclusion (\ref{e:inclusion}).  The idea is to use the Galois representations associated to Hilbert modular forms to construct an $R_p$-module $M$ and a surjection 
$\nabla_{\Sigma}^{\Sigma'}(H)_p^- \twoheadrightarrow M$ such that $\Fitt_{R_p}(M) \subset (\Theta)$.  
The properties of Fitting ideals imply that (\ref{e:inclusion}) follows from the existence of such a surjection.
As described in (\ref{e:nablasurj}),  a surjection from $\nabla_{\Sigma}^{\Sigma'}(H)_p^-$ can be constructed by defining a cohomology class $[\kappa] \in H^1(G_F, M)$ satisfying certain local conditions along with local trivializations at the places in $\Sigma$.

Ribet's method was described in great detail by Mazur in a beautiful article written for the celebration of Ribet's 60th birthday \cite{mazur}.   We borrow from this the following schematic diagram demonstrating the general path one follows to link $L$-values (in our case, the Stickelberger element $\Theta$) to class groups (in our case, the Ritter--Weiss module $\nabla_{\Sigma}^{\Sigma'}(H)_p^-$).
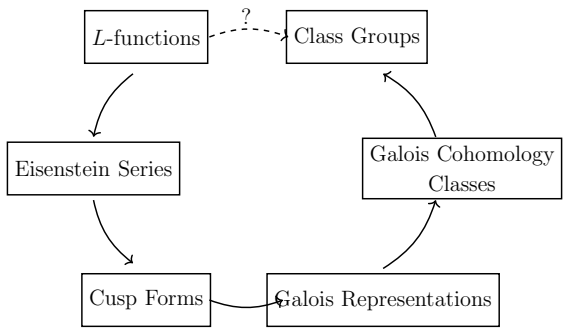
\begin{center}
\scalebox{0.7}{
\begin{tikzpicture}[every node/.style={minimum height={1cm},thick,align=center}]
\node[draw] (ES) at (-1, 0) {Eisenstein Series};
\node[draw] (LF) at (0, 2.5)  {$L$-functions};
\node[draw] (CG) at (4,2.5) {Class Groups};
\node[draw] (CF) at (0, -2.5){Cusp Forms};
\node[draw] (GC) at (6,0){Galois Cohomology \\ Classes};
\node[draw] (GR) at (4.5,-2.5) {Galois Representations};
\node at (1.9, 2.9) {?};
\draw[thick] (-0.25,1.8) edge [->, bend right =20]   (-1,0.6);
\draw[thick] (-1,-0.6) edge [->, bend right=20] (-0.25,-1.8);
\draw[thick] (1.2, -2.5) edge[->, bend right=20] (2.6,-2.5);
\draw[thick] (4.5,-1.9) edge[->, bend right=20] (5.5, -0.6);
\draw[thick] (5.5, 0.6) edge[->, bend right=20] (4.5,1.8);
\draw[thick,dashed] (1.2, 2.5) edge[->, bend left=20] (2.7,2.5);
\end{tikzpicture}}
\end{center}

Let us now trace this path in our application.

\subsubsection{$L$-functions to Eisenstein Series}  The key connection between $L$-functions and modular forms in Ribet's method is that $L$-functions appear as constant terms of Eisenstein series.  We  now define the space of modular forms in which our Stickelberger element $\Theta$ appears.

Let $k > 1$ be an integer such that \[ k \equiv 1 \pmod{(p-1)p^N} \] for a large value of $N$.  Let $\fn \subset O_F$ denote the conductor of $H/F$.  
  Let $M_k(\fn ; \Z)$ denote the group of Hilbert modular forms for $F$ of level $\fn $ with Fourier coefficients in $\Z$.  For each odd character $\chi$ of $G$ valued in $\C_p^*$, let \[ M_k(\fn , \chi) \subset M_k(\fn ; \Z) \otimes \C_p\] denote the subspace of forms of nebentypus $\chi$.
Let \[ \bchi\colon G_F \longrightarrow G \longrightarrow R_p^* \]
denote the canonical character, where the first arrow is projection and the second is induced by $G \hookrightarrow \Z[G]^*$.
\begin{definition} The space $M_k(\fn , \bchi; R_p)$ of {\em group-ring valued Hilbert modular forms of weight $k$ and level $\fn $ over $R_p$}  consists of those $f \in M_k(\fn ; \Z) \otimes R_p$ such that  
$\chi(f) \in M_k(\fn , \chi)$ for each odd character $\chi$.  Let $S_k(\fn , \bchi; R_p)$ denote the subspace of cusp forms. We define $M_k(\fn , \bchi; \Frac(R_p))$ and $S_k(\fn , \bchi; \Frac(R_p))$ similarly.  
\end{definition}

Hilbert modular forms $f$ are determined by their Fourier coefficients $c(\fm, f)$ indexed by the non-zero ideals $\fm \subset O_F$ and their constant terms $c_\lambda(0,f)$ indexed by $\lambda \in \Cl^+(F)$, the narrow class group of $F$.   For odd $k \ge 1$, there is an Eisenstein series $E_k(\bchi,1) \in M_k(\fn, \bchi)$ whose Fourier coefficients are given by
\[ c(\fm, E_k(\bchi, 1)) = \sum_{\substack{\fa \mid \fm \\(\fm/\fa, \fn)=1}} \!\!\!\! \bchi(\fm/\fa) \N\fa^{k-1}. 
\]
To describe the constant coefficients of $E_k(\bchi, 1)$ we first set $S$ to be minimal, i.e.\ the union of $S_\infty$ and $S_{\ram}$, where the latter is the set of primes dividing $\fn$.  Next we assume for the remainder of the article that $\fn \neq 1$; the case $\fn = 1$ causes no difficulties but the formulas must be slightly modified.  We then have
\[ c_\lambda(0, E_1(\bchi, 1)) =  \begin{cases} 0 & k > 1 \\ 2^{-n} \Theta^\#_{S} & k = 1. \end{cases} \]
Here $\Theta_S = \Theta_{S, \phi}(H/F, 0) \in \Q[G]$ denotes the $S$-depleted but unsmoothed Stickelberger element.
We have $E_k(\bchi, 1) \in M_k(\fn , \bchi; R_p)$ for $k > 1$ and $E_1(\bchi, 1) \in M_1(\fn , \bchi; \Frac(R_p))$ because of the possible non-integrality of the constant term.

\subsection{Eisenstein Series to Cusp Forms} In order to define a cusp form from the Eisenstein series, one is led to consider certain linear combinations of the analogues of $E_k(\bchi, 1)$ as $H$ ranges over all its CM subfields containing $F$.  This process also incorporates smoothing at the primes in $T$.  We avoid stating the slightly complicated formula here (see \cite{bsdk}*{Proposition 8.14}), but the end result is a group ring form $W_k(\bchi, 1)$ whose constant terms are given by
\begin{equation} \label{e:constant} c_\lambda(0, W_k(\bchi, 1)) = \begin{cases} 0 & k > 1 \\ 2^{-n} \Theta^\# & k = 1, \end{cases}
\end{equation}
where we remind the reader that $\Theta^\# = \Theta^\#_{\Sigma, \Sigma'}$.
 Building off the computations of \cite{acta}, we calculate in \cite{bsdk}*{\S8}  the constant terms of the $W_k(\bchi, 1)$ at {\em all cusps}; the  terms in (\ref{e:constant}) can be viewed as the constant terms ``at infinity.''  Indeed, it is the attempt to cancel the constant terms at other cusps that leads naturally to the definition of the $W_k(\bchi, 1)$.

In order to define a cusp form, we apply two important results of Silliman \cite{silliman}.  The first of these generalizes a result of Hida and Wiles and is stated below.
\begin{theorem}[\cite{silliman}*{Theorem 10.7}]  \label{t:hida} Let $m$ be a fixed positive integer. For positive integers $k \equiv 0 \pmod{(p-1)p^N}$ with $N$ sufficiently large, 
there is a Hilbert modular form $V_{k}$ of level 1, trivial nebentypus, and weight k defined over $\Z_p$ such that 
\[ V_{k} \equiv 1 \pmod{p^m}, \] and such that the normalized constant term  of $V_k$ at every cusp is congruent to $1 \pmod{p^m}$.
\end{theorem}

The idea to construct a cusp form is to fix a very large integer $m$ and to consider the product $W_1(\bchi, 1) V_k \in M_{k+1}(\fn, \bchi, R_p)$ with $V_k$ as in Theorem~\ref{t:hida}.  This series has constant terms at infinity  congruent to $2^{-n} \Theta^\#$ modulo $p^m$.  One then wants to subtract off $2^{-n} \Theta^\#H_{k+1}(\bchi)$  for some group ring valued form $H_{k+1}(\bchi) \in M_{k+1}(\fn, \bchi, R_p)$ to obtain a cusp form.  
If there exists a prime above $p$ dividing $\fn$ (i.e. $\Sigma - S_\infty$ is nonempty), then this strategy works.
Silliman's second result, which generalizes a result of Chai and is stated in \cite{silliman}*{Theorem 10.10}, implies that one can obtain a form that is cuspidal at the cusps ``above infinity at $p$'' in this fashion.  Applying Hida's ordinary operator then yields a form that is cuspidal.

\begin{theorem}[\cite{bsdk}*{Theorem 8.18}] \label{t:case2}
Suppose $\gcd(\fn, p) \neq 1$. For positive integers $k \equiv 1 \pmod{(p-1)p^N}$ and $N$ sufficiently large, there exists $H_k(\bchi) \in M_k(\fn, \bchi, R_p)$ such that
\[ {F}_k(\bchi) = e_p^{\ord}\left(W_1(\bchi, 1)V_{k-1} - \Theta^\#H_k(\bchi) \right) \]
lies in $S_k(\fn p, R, \bchi)$.
\end{theorem}
The significance of Theorem~\ref{t:case2} is that we have now constructed a cusp form that is congruent to an Eisenstein series modulo $\Theta^\#$.

When $\gcd(\fn, p)=1$, the construction of the cusp form is in fact more interesting, and a new feature appears.  In this case, the ordinary operator at $p$ does not annihilate the form $W_{k}(\bchi, 1)$, and it must be incorporated into our linear combination.  Moreover, this apparent cost has a great benefit---we obtain a congruence between a cusp form and Eisenstein series not only modulo $\Theta^\#$, but modulo a multiple $x \cdot \Theta^\#$ for a certain $x \in R_p$.  This element $x$ has an intuitive meaning---it represents the trivial zeroes of the $p$-adic $L$-function associated to $\bchi$, even the ``mod $p$ trivial zeroes.''  The precise definition is as follows.
\begin{lemma} \label{l:xdef} Suppose $\gcd(\fn, p)=1$.
For  positive $k \equiv 1 \pmod{(p-1)p^N}$ with  $N$ sufficiently large, the element
\[ x =  \frac{\Theta_{S_\infty}(1-k)}{ \Theta_{S_\infty}(0) } \in \Frac(R) \]
lies in $R$ and  is a non-zerodivisor.

\end{lemma}
The analogue of Theorem~\ref{t:case2} for $\gcd(\fn, p)=1$ is as follows.
\begin{theorem}[\cite{bsdk}*{Theorem 8.17}] \label{t:case1}
Suppose $\gcd(\fn, p) = 1$. For positive integers $k \equiv 1 \pmod{(p-1)p^N}$ and $N$ sufficiently large, there exists  $H_k(\bchi) \in M_k(\fn, \bchi, R_p)$ such that
\[ {F}_k(\bchi) = e_p^{\ord}\left( xW_1(\bchi, 1)V_{k-1} - W_k(\bchi, 1) -x \Theta^\#H_k(\bchi) \right) \]
lies in $S_k(\fn p, R, \bchi)$.
\end{theorem}

The extra factor of $x$ in our congruence between the cusp form $F_k(\bchi)$ and a linear combination of Eisenstein series  plays an extremely important role in showing that the Galois cohomology classes we construct are unramified at $p$.

We conclude this section by interpreting the congruences of Theorems~\ref{t:case2} and~\ref{t:case1} in terms of Hecke algebras.  We consider the Hecke algebra $\tilde{\T}$ generated over $R_p$ by the operators 
$T_\fq$ for primes $\fq \nmid \fn p$ and $U_\fp$ for primes $\fp \mid p$.  (We ignore the operators $U_\fq$ for $\fq \mid \fn, \fq \nmid p$ in order to avoid issues regarding non-reducedness of Hecke algebras arising from the presence of {\em oldforms}.) We denote by $\T = e_p^{\ord}(\tilde{\T})$ Hida's ordinary Hecke algebra associated to $\tilde{\T}$.  
Let $\epsilon_{\cyc}\colon G_F \longrightarrow \Z_p^*$ denote the $p$-adic cyclotomic character of $F$.
Theorems~\ref{t:case2} and~\ref{t:case1} then yield:
\begin{theorem} \label{t:yxt}
Let $x = 1$ if $\gcd(\fn,p) \neq 1$ and let $x$ be as in Lemma~\ref{l:xdef} if $\gcd(\fn, p)=1$.
 There exists an $R_p/x\Theta^\#$-algebra $W$ and a surjective $R_p$-algebra homomorphism $\varphi\colon {\T} \longrightarrow W$ satisfying the following properties:
\begin{itemize}
\item The structure map $R_p/x\Theta^\# \longrightarrow W$ is an injection.
\item $\varphi(T_\fl) = \epsilon_{\cyc}^{k-1}(\fl) + \bchi(\fl)$  for  $\fl \nmid \fn p.$
\item $\varphi(U_\fp) = 1$  for  $\fp \mid \gcd(\fn, p).$ 
\item Let \[ {U} = \prod_{\fp \mid p, \ \fp \nmid \fn} (U_{\fp} - \bchi(\fp)) \in {\T}. \]
  If $y \in R_p$ and $\varphi(U)y = 0$ in $W$, then $y \in (\Theta^\#)$. 
\end{itemize}
\end{theorem}
The idea of this theorem is the usual one: the homomorphism $\varphi$ sends a Hecke operator to its ``eigenvalue mod $x \Theta^\#$'' acting on $F_k(\bchi)$.  The only subtlety is the last bullet point: the operators $U_\fp$ for $\fp \mid p, \  \fp \nmid \fn$ do not act as scalars, so a more involved argument is necessary. This explains why the ring $W$ is not just $R_p/x \Theta^\#$.  
The idea behind the last statement of the theorem is that the operator $\varphi(U)$ introduces a factor of $x$; If $xy$ is divisible by $x\Theta^\#$ in $R_p$, then 
$y$ is divisible by $\Theta^\#$ since $x$ is a non-zerodivisor.  This demonstrates the essential additional ingredient provided by the ``higher congruence'' modulo $x \Theta^\#$ rather than just modulo $\Theta^\#$.
See \cite{bsdk}*{Theorem 8.23} for details.

\subsection{Cusp Forms to Galois Representations}

In this section we study the Galois representation attached to cusp forms that are congruent to Eisenstein series.  Let $\fm$ be the intersection of the finitely many maximal ideals of $\T$ containing the kernel of $\varphi$.  Put $\T_{\fm}$ for the completion of $\T$ with respect to $\fm$ and $K = \Frac(\T_{\fm})$. Then $K$ is a finite product of fields parameterized by the $\Q_p$-Galois orbits of cuspidal newforms of weight $k$ and level $\fn$, defined over the ring of integers in a finite extension of $\Z_p$, that are congruent to an Eisenstein series modulo the maximal ideal.
As in \cite{bsdk}*{\S9.2}, the work of Hida and Wiles gives a Galois representation
\[
\rho : G_F \rightarrow \GL_2(K)
\]
satsifying the following.  \begin{enumerate}
\item $\rho$ is unramified outside $\fn p$. 
\item For all primes $\fl \nmid \fn p$, the characteristic polynomial of $\rho(\Frob_{\fl})$ is given by 
\begin{equation} \label{e:charpoly}
\text{char}(\rho(\Frob_{\fl}))(x) = x^2 - T_{\fl} x + \bchi(\fl) \N \fl^{k-1}.
\end{equation}
\item For every $\fq \mid p$, let $G_{\fq}$ denote a decomposition group at $\fq$. We have
\[
\rho|_{G_{\fq}} \sim \left( \begin{array}{cc} \bchi \varepsilon_{\cyc}^{k-1} \eta_{\fq}^{-1} & \ \  * \\ 0 & \ \ \eta_{\fq} \end{array} \right),
\]
where $\varepsilon_{\cyc}$ is the $p$-adic cyclotomic character and $\eta_{\fq}$ is an unramified character given by $\eta_{\fq}(\rec_\fq(\varpi^{-1})) = U_{\fq}$, with $\varpi$ a uniformizer of $F_\fq^*$. 
\end{enumerate}
Let $\I$ denote the kernel of $\varphi$ extended to $\varphi: \T_{\fm} \rightarrow W$.  
Reducing (\ref{e:charpoly}) modulo $\I$, and using \v{C}ebotarev to extend from $\Frob_{\fl}$ to all $\sigma \in G_F$, 
we see that the characteristic polynomial of $\rho(\sigma)$ is congruent to $(x - \bchi(\sigma))(x- \epsilon_{\cyc}(\sigma)) \pmod{\I}$.
In particular if $\chi(\sigma) \not\equiv \epsilon_{\cyc}(\sigma) \pmod{\fm}$, then by Hensel's Lemma $\rho(\sigma)$ has distinct eigenvalues $\lambda_1, \lambda_2 \in \T_{\fm}$ such that $\lambda_1 \equiv \epsilon_{\cyc}^{k-1}(\sigma)  \pmod{I}$ and $\lambda_2 \equiv \bchi(\sigma)  \pmod{\I}$.

To define a convenient basis for $\rho$, we choose $\tau \in G_F$ such that:
\begin{enumerate}
\item $\tau$ restricts to the complex conjugation of $G$, 
\item for each $\fq \mid p$, the eigenspace of $\rho|_{G_{\fq}}$ projected to each factor of $K$ is not stable under $\rho(\tau)$.
\end{enumerate} 
See \cite{bsdk}*{Proposition 9.3} for the existence of such $\tau$.  Since $p \neq 2$, we have \[ \bchi(\tau) = - 1 \not\equiv 1 \equiv \epsilon_{\cyc}(\tau) \pmod{\fm}. \]
It follows from the discussion above that the eigenvalues of $\rho(\tau)$ satisfy 
 $\lambda_1 \equiv \epsilon_{\cyc}(\sigma) \pmod{\I}$ and $\lambda_2 \equiv -1 \pmod{\I}$.
Fix the basis consisting of eigenvectors of $\rho(\tau)$, say $\rho(\tau) =\begin{pmatrix}
 \lambda_1 & \ \ 0 \\ 0 & \ \  \lambda_2\end{pmatrix}$.
  For a general $\sigma \in G_F$, write 
$\rho(\sigma) = \begin{pmatrix} a(\sigma) & \ \ \ b(\sigma) \\ c(\sigma) &  \ \ \ d(\sigma) \end{pmatrix}.$ 
For each $\fq \mid p$, there is a change of basis matrix $M_{\fq} = \begin{pmatrix}
 A_{\fq} & \ \  B_{\fq} \\ C_{\fq} & \ \ D_{\fq} \end{pmatrix} \in \GL_2(K)$ satisfying
\begin{equation} \label{e:changebasis}
\begin{pmatrix}
a(\sigma) & \ \ \ b(\sigma) \\ c(\sigma) & \ \ \ d(\sigma) \end{pmatrix}
 M_{\fq} = M_{\fq} \begin{pmatrix}
  \bchi \varepsilon_{\cyc}^{k-1} \eta_{\fq}^{-1} & \ \ * \\ 0 & \ \ \eta_{\fq} \end{pmatrix}.
\end{equation}
The second condition in the choice of $\tau$ ensures that $A_{\fq}, C_{\fq} \in K^*$. Furthermore, equating the upper left hand entries in (\ref{e:changebasis}) gives:
\begin{equation} \label{e:brelation}
b(\sigma) = \frac{A_{\fq}}{C_{\fq}} ( a(\sigma) - \bchi \varepsilon_{\cyc}^{k-1} \eta_{\fq}^{-1} (\sigma)) \qquad \quad \text{ for all } \sigma \in G_{\fq}.
\end{equation}

\subsection{Galois Representations to Galois Cohomology Classes}
 
We summarize \cite{bsdk}*{\S 9.3}.
As explained above, we have \begin{equation} \label{e:adeq}
 a(\sigma) + d(\sigma) \equiv \bchi(\sigma) + \epsilon_{\cyc}^{k-1}(\sigma) \pmod{I} \quad \text{ for all } \sigma \in G_F. 
 \end{equation}
  Applying the same rule for $\tau \sigma$ and noting that $a(\tau \sigma) = \lambda_1 a(\sigma), \  d(\sigma \tau) = \lambda_2 d(\sigma)$, we find
\begin{equation} \label{e:adeq2}
 a(\sigma)\epsilon^{k-1}_{\cyc}(\tau) -  d(\sigma) \equiv - \bchi(\sigma) + \epsilon_{\cyc}^{k-1}(\sigma\tau) \pmod{I}. 
 \end{equation}
Solving the congruences (\ref{e:adeq}) and (\ref{e:adeq2}) and once again using the fact that $\epsilon_{\cyc}^{k-1}(\tau) \not\equiv -1 \pmod{\fm}$ since $p \neq 2$, we find
 that $a(\sigma), d(\sigma) \in \T_{\fm}$ and 
\begin{equation} \label{e:adcong}
a(\sigma) \equiv \varepsilon_{\cyc}^{k-1}(\sigma) \pmod{\I}, \quad  d(\sigma) \equiv \bchi(\sigma) \pmod{\I} \quad
\text{ for all }\sigma \in G_F.
\end{equation}

Let $B$ be the $\T_{\fm}$ submodule of $K$ generated by $\{b(\sigma): \sigma \in G_F\} \cup \{\frac{A_{\fq}}{C_{\fq}} : \fq \in \Sigma \setminus S_{\infty}\}$. We have $\rho(\sigma \sigma') = \rho(\sigma) \rho(\sigma')$ for $\sigma, \sigma' \in G_F$. Equating the upper right entries and using equation (\ref{e:adcong}), we obtain
\begin{equation} \label{e:bcocycle}
b(\sigma \sigma') = a(\sigma) b(\sigma') + b(\sigma) d(\sigma') \equiv \varepsilon_{\cyc}^{k-1}(\sigma) b(\sigma') + \bchi(\sigma') b(\sigma) \pmod{\I B}. 
\end{equation}
Let $m$ be am integer such that $k \equiv 1 \pmod{(p-1)p^m}$. Let $I_{\fq}$ denote the inertia subgroup of $G_F$ of a prime $\fq$. Put $B_1$ for the $\T_\fm$-submodule of $B$ generated by 
\[
\I B \cup p^m B \cup  \{b(\sigma) : \sigma \in I_\fq \text{ for } \fq \mid p, \ \fq \notin \Sigma \}.
\] 
Define $\overline{B} = B/B_1$. Equation (\ref{e:bcocycle}) then gives that $\kappa(\sigma) = \bchi^{-1}(\sigma) b(\sigma)$ is a cocyle defining a cohomology class $[\kappa]$ in $H^1(G_F, \overline{B}(\bchi^{-1}))$ satisfying the following local properties.
\begin{enumerate}
\item As $\rho$ is unramified at $\fl \nmid \fn p$, so is the class $[\kappa]$. 
\item As $\overline{B}$ is pro-$p$, the class $[\kappa]$ is at most tamely ramified at any prime $\fl \mid \fn$ not above $p$.  
\item It is proven in \cite{bsdk}*{\S 4.1} that we may assume $\Sigma'$ does not contain any primes above $p$. Thus $[\kappa]$ is at most tamely ramified at all primes in $\Sigma'$.
\item By the definition of $B_1$, where we have included $b(I_\fq)$ for primes $\fq \mid p, \ \fq \not\in \Sigma$, the class $[\kappa]$ is unramified at such $\fq$.
\item Equation (\ref{e:brelation}) implies that $[\kappa]$ is locally trivial at finite primes in $\Sigma$. As $p$ is odd, $[\kappa]$ is locally trivial at archimedian places (\cite{bsdk}*{Proposition 9.5}).
\end{enumerate}

\subsection{Galois Cohomology Classes to Class Groups}

The Galois cohomology class $[\kappa]$ satisfies the conditions listed after equation (\ref{e:nablasurj}) and hence gives a surjection 
\[
\nabla^{\Sigma'}_{\Sigma}(H)^{-}_p \longtwoheadrightarrow \overline{B}(\bchi^{-1}).
\]
For details see \cite{hdk}*{Theorem 4.4}. The general properties of Fitting ideals imply 
\[
\Fitt_{R_p}(\nabla^{\Sigma'}_{\Sigma}(H)^{-}_p) \subset \Fitt_{R_p}(\overline{B}(\bchi^{-1})).
\]
It is therefore enough to prove that $\Fitt_{R_p}(\overline{B}) \subset (\Theta^{\#})$. Typically in Ribet's method, one argues that the fractional ideal $B$ is a faithful $\T_\fm$-module, and hence the Fitting ideal of $B/ \I B$ is contained in $\I$. However, our module 
$\overline{B}$ is more complicated than $B/\I B$, so we proceed as follows. Using equation (\ref{e:brelation}), we show that any element in $\Fitt_{R_p}(\overline{B})$ is annihilated by $\varphi(U)$ for the operator $U$ from Theorem \ref{t:yxt}. The final assertion of this theorem then implies that $\Fitt_{R_p}(\overline{B})$ contained in $(\Theta^\#)$. See \cite{bsdk}*{\S9.5} for details.

This concludes our summary of the proof of Theorem~\ref{t:selfit}.

\section{Explicit Formula for Brumer--Stark Units}

In this final section of the paper, we discuss the first author's explicit formula for Brumer--Stark units as mentioned in \S\ref{s:ef}.
The conjecture in the case that $F$ is a real quadratic field was studied in \cite{darmondasgupta}, and the general case was studied in \cite{dasguptaduke}.  Here we consider an arbitrary totally real field $F$, but to simplify formulas we assume that the rational prime $p$ is inert in $F$.  Furthermore, we let $H$ be the narrow ray class field of some conductor $\fn \subset O_F$ and assume that $p \equiv 1\pmod{\fn}$.  This ensures that the prime  $\fp = p O_F$ splits completely in $H$.  Fix a prime $\fP$ of $H$ above $\fp$.  
We fix $S \supset S_\infty \cup S_\ram = \{v \mid \fn \infty\}$.
We also fix a prime ideal $\fl \subset O_{F}$ such that $\N\fl = \ell > n+1$ is a prime integer and let $T = \{\fl\}$.

In this setting, we will present a $\fp$-adic analytic formula for the image of the Brumer--Stark unit $u_\fp \in H^*$ in $H_\fP^* \cong F_\fp^*$.
The most general, conceptually satsifying, and theoretically useful form of this conjecture uses the Eisenstein cocycle.  This is a class in the $(n-1)$st cohomology of $\GL_n(\Z)$ that has many avatars studied by several authors (see \cites{bcg, charolloisdasgupta, cdg, darmondasgupta, dasguptaspiess, bkl, spiess}).
In this paper, we avoid defining the Eisenstein cocycle and present instead the more explicit and down to earth version of the conjectural formula for $u_\fp$ stated in \cite{dasguptaduke}.

\subsection{Shintani's Method}
 Fixing an ordering of the $n$ real embeddings of $F$ yields an map $F \hookrightarrow \R^n$ such that the image of any fractional ideal is a cocompact lattice.  We let $F^*$ act on $\R^n$ by composing this embedding with componentwise multiplication and denote the action by $*$.

Let $v_1, \dotsc, v_r \in (\R^{>0})^n$, $1 \le r\le n$, be vectors in the totally positive orthant that are linearly independent over $\R$.  The corresponding {\em simplicial cone} is defined by
\[ C(v_1, \dotsc, v_r) = \left\{ \sum_{i=1}^{r} t_i v_i \colon 0 < t_i \right\} \subset (\R^{>0})^n. \]
Suppose now $r=n$.
We will define a certain union of $C(v_1, \dotsc, v_n)$ and some of its boundary faces that we call the {\em Colmez closure}.
Write \[ (0, 0, \dotsc, 1) = \sum_{i=1}^{n} q_i v_i, \qquad q_i \in \R. \]
For each nonempty subset $J \subset \{1, \dotsc, n\},$ we say that $J$ is {\em positive} if $q_i > 0$ for all $i \not \in J$.  The Colmez closure of $C(v_1, \dotsc, v_n)$ is defined by:
\[ C^*(v_1, \dotsc, v_n) = \bigsqcup_{J \text{ positive}} C(\{v_j, j \in J\}). \]

Let $E(\fn) \subset O_F^*$ denote the subgroup of totally positive units $\epsilon$ such that $\epsilon \equiv 1 \pmod{\fn}$.  
Shintani proved that there exists a union of simplicial cones that is a fundamental domain for the action of $E(\fn)$ on 
$(\R^{>0})^n$.  For example, in the real quadratic case ($n=2$), $E(\fn) = \langle \epsilon \rangle$ is cyclic and $C^*(1, \epsilon)$ is a fundamental domain.  In the general case, it can be difficult to write down an explicit fundamental domain, but a nice generalization of the $n=2$ case is obtained if we allow ourselves to consider instead a {\em signed fundamental domain}.  For a simplicial cone $C$, let ${\mathbf 1}_C$ denote the characteristic function of $C$ on $(\R^{>0})^n$.
\begin{definition} \label{d:sfd}
A {\em signed fundamental domain}
for the action of $E(\fn)$ on $(\R^{>0})^n$ is by definition a formal linear combination $D = \sum_i a_i C_i$ of
simplicial cones $C_i$ with $a_i \in \Z$ such that \[ \sum_{u \in E(\fn)} \sum_i a_i {\mathbf 1}_{C_i}(u *x) = 1 \]
for all $x \in (\R_{>0})^n$.  
\end{definition}

Fix an ordered basis $\{\epsilon_1, \dotsc, \epsilon_{n-1}\}$ for $E(\fn)$.
Define the orientation
\begin{equation} \label{e:wedef}
 w_\epsilon = \sign\det( \log(\epsilon_{ij}))_{i,j = 1}^{n-1}) = \pm 1, 
 \end{equation}
where $\epsilon_{ij}$ denotes the $j$th coordinate of  $\epsilon_i$.
 For each permutation $\sigma \in S_{n-1}$ let
\[ v_{i, \sigma} = \epsilon_{\sigma(1)} \cdots \epsilon_{\sigma(i-1)} \in E(\fn), \qquad i = 1, \dotsc, n. \]
By convention $v_{1, \sigma} = (1, 1, \dotsc, 1)$ for all $\sigma$.
Define
\[ w_\sigma = (-1)^{n-1} w_\epsilon \sign(\sigma) \sign(\det(v_{i,\sigma})_{i=1}^{n}) \in \{0, \pm 1\}. \]

The following result was proven independently by Diaz y Diaz--Friedman \cite{ddf} and Charollois--Dasgupta--Greenberg \cite{cdg}*{Theorem 1.5}, generalizing the result of Colmez \cite{colmez} in the case that all $w_\sigma = 1$.
\begin{theorem}   \label{t:ddf}
The formal linear combination
\[ \sum_{\sigma \in S_{n-1}} w_\sigma C^*(v_{1, \sigma}, \dotsc, v_{n, \sigma})\]
is a signed fundamental domain for the action of $E(\fn)$ on $(\R^{>0})^n$.
\end{theorem}

\subsection{The Formula} \label{s:exact}

Throughout this section assume that $p$ is odd.  Recall that $T = \{\fl\}$.
Let $\fb$ be a fractional ideal that is relatively prime to $\fn \fl$.  Let $D = \sum a_i C_i$ be the signed fundamental domain for the action of $E(\fn)$ on $(\R^{>0})^n$ given in Theorem~\ref{t:ddf}.  We use all this data to define a $\Z$-valued measure $\mu$ on $O_p$, the $p$-adic completion of $O_F$.
Fix an element $z \in \fb^{-1}$ such that $z \equiv 1 \pmod{\fn}$.
For each compact open set $U \subset O_p$, define the {\em Shintani zeta function}
\[ \zeta(\fb, U, D, s) = \sum_i a_i \!\!\!\! \sum_{\substack{\alpha \in C_i \cap \fb^{-1}\fn + z \\ \alpha \in U, (\alpha, S)=1}} (\N\alpha)^{-s}. \]
Shintani proved that this sum converges for $\Re(s)$ large enough and extends to a meromorphic function on $\C$.  Define
\[ \mu_\fb(U) = \zeta(\fb, U, D, 0) - \ell \cdot \zeta(\fb\fl^{-1}, U, D, 0).
\]
Using Shintani's formulas, one may show:
\begin{theorem}[Proposition 3.12, \cite{dasguptaduke}]  For every compact open $U \subset O_p$, we have $\mu_{\fb}(U) \in \Z$.
\end{theorem}

We may now state our conjectural exact formula for the Brumer--Stark unit $u_p$ and all of its conjugates over $F$.
Write \[ \Theta_{S,T} = \sum_{\sigma \in G} \zeta_{S,T}(\sigma) \sigma^{-1}, \quad \zeta_{S,T}(\sigma) \in \Z. \]
Define
\begin{equation} \label{e:explicit}
u_{p}(\fb)^{\an} = p^{\zeta_{S,T}(\sigma_\fb)} \mint_{O_p^*} x \ d\mu_\fb(x) \in F_p^*. 
\end{equation}
Here the crossed integral is a {\em multiplicative integral} in the sense of Darmon \cite{darmon} and can be expressed as a limit of Riemann products:
\[  \mint_{O_p^*} x \ d\mu_\fb(x)  := \lim_{m \rightarrow \infty} \prod_{a \in (O_p/p^m)^*} a^{\mu_\fb(a + p^mO_p)}. 
\]
Write $\sigma_\fb \in G$ for the Frobenius associated to $\fb$.  In \cite{dasguptaduke}*{Theorem 5.15} we prove that $u_p(\fb)^{\an}$ depends only on the image of $\fb$ in the narrow ray class group of conductor $\fn$, i.e.\ on $\sigma_\fb \in G$ (at least up to a root of unity in $F_p^*$).
\begin{conjecture}  \label{c:explicit}  
 We have
$
 \sigma_\fb(u_\fp) = u_{p}(\fb)^{\an}$ in $F_p^*$.
 \end{conjecture}
 The expression (\ref{e:explicit}) can be computed to high $p$-adic precision on a computer.  See~\cite{fl} for tables of narrow Hilbert class fields of real quadratic fields determined using this formula. 
 
It is convenient to have an invariant that also satisfies $u_p(\fb \fq)^{\an} = (u_p(\fb)^{\an})^{-1}$ if $\fq$ is a prime such that $\sigma_\fq = c$.  Conjecture~\ref{c:explicit} would imply such a formula, but it is unclear whether this purely analytic statement can be proved unconditionally.
To this end we fix $\fq$ such that $\sigma_\fq = c$ and define 
\[ v_{p}(\fb)^{\an} = \left( \frac{u_{p}(\fb)^{\an}}{u_{p}(\fb \fq)^{\an}}\right)^{1/2} \in \hat{F}_p^* := F_p^* \hat{\otimes} \Z_p. \]
One then has \begin{equation} \label{e:vinv}
v_p(\fb \fq)^{\an} = (v_p(\fb)^{\an})^{-1}\end{equation} unconditionally, and we expect to have $v_p(\fb)^{\an} = u_p(\fb)^{\an}.$  The following is therefore a slightly easier form of Conjecture~\ref{c:explicit} to study.
\begin{conjecture}  \label{c:explicit2}  
 We have
$
 \sigma_\fb(u_\fp) = v_{p}(\fb)^{\an}$ in $\hat{F}_p^*$.
 \end{conjecture}

\subsection{Horizontal Iwasawa Theory} \label{s:hit}

We now discuss the relationship between Gross's tower of fields conjecture (Conjecture~\ref{c:grosstower}) and our conjectural exact formula for Brumer--Stark units. 
Our goal is to prove:
\begin{theorem} \label{t:gimplyh}
Assume that $p$ is odd.  Gross's conjecture implies Conjecture~\ref{c:explicit2}.
\end{theorem}
 
 In this exposition we have assumed that the odd prime $p$ is inert in $F$ and that $\fp = p O_F$.  In the case of general $\fp$, one must still assume that $p$ is odd and unramified in $F$ in the statement of Theorem~\ref{t:gimplyh}.
 
The abelian extensions $L/F$ to which we can apply Gross's conjecture (with $S' = S \cup \{\fp\})$ as in Conjecture~\ref{c:grosstower}) are those that contain $H$ and are unramified outside $S' \infty$.  Let $F_{S'}$ denote the maximal abelian extension of $F$ unramified outside $S' \infty$.  The reciprocity map of class field theory yields an explicit description of 
$\Gal(F_{S'}/H)$.  For each finite $v \in S'$, let $U_{v, \fn} \subset O_v^*$ denote the subgroup of elements congruent to 1 modulo $\fn O_v$ (so $U_{f, \fn} = O_v^*$ for $v \nmid \fn$).  Define $\bO^* = \prod_{v \in S' \setminus S_\infty} U_{v, \fn}$.
Then \[ \Gal(F_{S'}/H) \cong \bO^* / \overline{E(\fn)}, \] where $\overline{E(\fn)}$ denotes the topological closure of ${E(\fn)}$ embedded diagonally in  $\bO^*$.

 For each finite extension $L \subset F_{S'}$ containing $H$, if we write $\Gamma = \Gal(L/H)$, then (\ref{e:recg}) yields a formula for 
 $\rec_G(u_p)$ in $I/I^2 \cong \Z[G] \otimes \Gamma. $
 Under this isomorphism, the coefficient of $\sigma_{\fb}^{-1}$ is just the image of $\rec_{\fp}(\sigma_\fb(u_p))$ in $\Gamma$.
 Taking the inverse limit over all $H \subset L \subset F_{S'}$ therefore gives an equality for 
 \[  (\sigma_\fb(u_p),1, 1, \dots, 1) \text{ in } \bO / \overline{E(\fn)}. \]
Here we have written $O_p^*$ as the first component of $\bO$.

The next key point is that the constructions of 
\S\ref{s:exact} can be repeated to provide a measure $\mu_{\fb, \bO}$ on $\bO =  \prod_{v \in S'\setminus S_\infty} O_{v}$ extending the measure $\mu_\fb$ on $O_p$. 
It is not hard to check that the restriction of $\mu_{\fb, \bO}$ to $\bO^*$, pushed forward to $\bO^*/\overline{E(\fn)}$, is precisely the measure that recovers the values of the partial zeta functions of the abelian extensions $L$ contained in $F_{S'}$.  These are exactly the values appearing in Gross's conjecture.  In other words, Gross's conjecture for the set $S'$ is equivalent to
\begin{equation} \label{e:shintaniS}
 (\sigma_\fb(u_\fp),1, 1, \dots, 1) \cdot p^{-\zeta_{S,T}(\sigma_\fb)}   =  \mint_{\bO^*} x \ d\mu_{\fb,\bO}(x) \quad \text{ in } \bO / \overline{E(\fn)}.
\end{equation}
See \cite{dasguptaduke}*{Proposition 3.4}.
The next important calculation (\cite{dasguptaduke}*{Theorem 3.22}) is that
\begin{equation} \label{e:mucompat}
p^{\zeta_{S,T}(\sigma_\fb)} \mint_{\bO^*} x \ d\mu_{\fb,\bO}(x)  = (u_p(\fb)^{\an}, 1, 1, \dots, 1).
\end{equation}
The first component of this is simply the compatibility of the constructions of $\mu_{\fb}$ and $\mu_{\fb, \bO}$; the interesting part of the computation is the 1's in the components away from $p$.  Equations (\ref{e:shintaniS}) and (\ref{e:mucompat}) combine to yield that the ratio $\sigma_\fb(u_p)/u_p(\fb)^{\an}$ lies in the group
\[ D(S) = \{x \in O_p^* \colon (x, 1, 1, \dotsc, 1) \in \overline{E(\fn)} \subset \bO^* \}. \]
We can also conclude  \begin{equation} \label{e:uvds}
\sigma_\fb(u_p)/v_p(\fb)^{\an} \in D(S) \end{equation}
 since $c(u_p) = u_{p}^{-1}$.

The final trick, inspired by the method of Taylor--Wiles, is to consider certain  enlarged sets $S_Q = S \cup Q$ for a well-chosen finite set of auxiliary primes $Q$. 
 Let us compare the Brumer--Stark units for $S$ and $S_Q$, denoted $u_{p}$ and $u_{p}( S_Q)$, respectively.  The defining property (\ref{e:bs}) shows that
\[ u_{p}(S_Q) = u_{p}^z, \quad \text{ where } z = \prod_{\fq \in Q}(1 - \sigma_{\fq}^{-1}) \in \Z[G]. \]
In particular if we choose the $\fq \in Q$ such that $\sigma_\fq = c$, the complex conjugation of $G$, then $u_{p}(S_Q) = u_{p}^{2^{\#Q}}.$  Using (\ref{e:vinv}), one can similarly show that 
$v_{p}(S_Q, \fb) = v_p(\fb)^{2^{\#Q}}$.
Now (\ref{e:uvds}) for $S_Q$ implies that \[ \sigma_\fb(u_p(S_Q))/v_p(S_Q, \fb)^{\an} \in D(S_Q), \]
hence \[ \left(\sigma_\fb(u_p) /v_p(\fb)^{\an}\right)^{2^{\#Q}} \in D(S_Q), \quad \text{ so } \quad
\sigma_\fb(u_p)/v_p(\fb)^{\an} \in D(S_Q)  \]
since $p \neq 2$.

To conclude the proof of~Theorem~\ref{t:gimplyh},
 one shows using the \v{C}ebotarev Density Theorem that one can choose the sets $Q$ to force $D(S_Q)$ as small as desired (i.e.\ the intersection of $D(S_Q)$ over all possible $Q$ is trival).    See~\cite{dasguptaduke}*{Lemma 5.17} for details.
 
\subsection{The Greenberg--Stevens $\sL$-invariant} \label{s:gs}

We briefly summarize our proof of the $p$-part of Gross's conjecture (Theorem~\ref{t:maingross}), which as just explained implies our explicit formula for Brumer--Stark units given in Conjecture~\ref{c:explicit2}.

The work of Greenberg and Stevens \cite{gs} was a seminal breakthrough in the study of trivial zeroes of $p$-adic $L$-functions.
Their perspective was highly influential in \cite{ddp}, where the rank one $p$-adic Gross--Stark conjecture was interpreted as the equality of an algebraic $L$-invariant $\sL_{\alg}$ and
an  analytic $L$-invariant $\sL_{\an}$.  The analytic $\sL$-invariant is the ratio of the leading term of the $p$-adic $L$-function at $s=0$ to its classical counterpart: \begin{equation} \label{e:lan}
 \sL_{\an} = -\frac{L_p'(\chi\omega, 0)}{L(\chi, 0)}. \end{equation}
   The algebraic $L$-invariant is the ratio of the $p$-adic logarithm and valuation of the $\chi^{-1}$-component of the Brumer--Stark unit: 
   \begin{equation} \label{e:lalg}
    \sL_{\alg} = \frac{\log_p \Norm_{H_{\fP}/\Q_p}(u_\fp^{\chi^{-1}}) }{ \ord_\fP(u_\fp^{\chi^{-1}})}.  
    \end{equation}

There is no difficulty in defining the ratios (\ref{e:lan}) and (\ref{e:lalg}), since the quantities live in a $p$-adic field and the denominators are non-zero.  The analogue of this situation for Gross's Conjecture~\ref{c:grosstower} is more delicate. The role of the $p$-adic $L$-function is played by the Stickelberger element 
$\Theta_L := \Theta_{S',T}(L/F,0) \in \Z[\fg]$, and the analogue of the derivative at 0 is played by the image of $\Theta_L$ in $I/I^2$.  The role of the classical $L$-function is played by the element $\Theta_H := \Theta_{S, T}(H/F,0) \in \Z[G]$.  It is therefore not clear how to take the ``ratio'' of these quantities.  Similarly, the role of the $p$-adic logarithm is played by $\rec_G(u_{\fp}) \in I/I^2$ and the role of the $p$-adic valuation is played by $\ord_G(u_{\fp}) \in \Z[G]$.

For this reason, we introduce in \cite{hdk} an $R$-algebra $R_\sL$ that is generated by an element $\sL$ that plays the role of the analytic 
$\sL$-invariant, i.e. the ``ratio" between 
$\Theta_L$ and $\Theta_H$.
We define
\begin{equation} \label{e:rldef}
 R_\sL = R[\sL] / (\Theta_H \sL - \Theta_L, \sL I, \sL^2, I^2). 
 \end{equation}
 A key nontrivial result is that this ring, in which we have adjoined a ratio $\sL$ between $\Theta_L$ and $\Theta_H$, is still large enough to see  $R/I^2$.  
 
 \begin{theorem}[\cite{hdk}, Theorem 3.4] The kernel of the structure map $R \longrightarrow R_{\sL}$ is $I^2$. \label{t:ws}
\end{theorem}

It follows from this theorem that Gross's Conjecture is equivalent to the equality \begin{equation} 
 \rec_G(u_{\fp})  = \sL \ord_G(u_{\fp}) \text{ in } R_{\sL}, \label{e:linv} \end{equation}
since the right side is by definition $\sL \Theta_H = \Theta_L$.

 To prove (\ref{e:linv}), we define a generalized Ritter--Weiss module $\nabla_{\!\sL}$ over the ring $R_{\sL}$  that can be viewed as a gluing of the modules $\nabla_{S}^{T}(H)$ and $\nabla_{S'}^{T}(L)$.
 We  show in \cite{hdk}*{Theorem 4.6} that the Fitting ideal $\Fitt_{R_\sL}(\nabla_{\!\sL})$ is generated by the element 
\[ \rec_G(u_{\fp}) -\sL \ord_G(u_{\fp}) \in I/I^2,  \]
and hence that (\ref{e:linv}) is equivalent to \begin{equation} \label{e:ifit}
 \Fitt_{R_\sL}(\nabla_{\!\sL}) = 0. 
  \end{equation}
  (For the sake of accuracy, we remark that in reality we do all of this with  $(S,T)$ replaced by the pair $(\Sigma, \Sigma')$ defined in \S\ref{s:burns}, as in \S\ref{s:rw}.)
  
  The vanishing of $ \Fitt_{R_\sL}(\nabla_{\!\sL})$ is proven following the methods of \S\ref{s:rw}.  
  We interpret surjective homomorphisms from $\nabla_{\sL}$ to $R_{\sL}$-modules $M$ in terms of Galois cohomology classes satisfying certain local conditions.  We construct a suitable  Galois cohomology class valued in a module $M$ using 
  an explicit construction with group-ring valued Hilbert modular forms and their associated Galois representations.  The module $M$ is shown to be large enough that its Fitting ideal over $R_{\sL}$ vanishes, whence the same is true for $\nabla_{\sL}$ since it has $M$ as a quotient.

\subsection{The method of Darmon--Pozzi--Vonk} \label{s:dpv}

We conclude by describing a proof of Conjecture~\ref{c:explicit} in the case that $F$ is a real quadratic field
 in the beautiful work of Darmon, Pozzi, and Vonk \cite{dpv}.  Their method is purely $p$-adic (i.e. ``vertical''), rather than involving the introduction of auxiliary primes (i.e. ``horizontal'').  The strategy follows a rich history of arithmetic formulas proven by exhibiting both sides of an equation as certain Fourier coefficients in an equality of modular forms.  For instance, Katz gave an elegant proof of Leopoldt's evaluation of the Kubota--Leopoldt $p$-adic $L$-function at $s=1$ by exhibiting an equality of $p$-adic modular forms, one of whose constant terms is the $p$-adic $L$-value and the other is the $p$-adic logarithm of a unit (see \cite{katz}*{\S10.2}).  The proof of Darmon--Pozzi--Vonk follows a similar strategy.

Let $F$ be a real quadratic field, $p$ an odd prime, and $H$ a narrow ring class field extension of $F$ (so in particular $p O_F$ splits completely in $H$).
Darmon--Pozzi--Vonk demonstrate an  equality of certain classical modular forms of weight 2 on $\Gamma_0(p) \subset \SL_2(\Z)$ that we denote $f_1$ and $f_2$.

This first of these forms $f_1$ is obtained by considering a Hida family of Hilbert modular cusp forms for $F$ specializing in weight 1 to a $p$-stabilized Eisenstein series.  The constant term of this weight 1 Eisenstein series vanishes because of the trivial zero of the corresponding $p$-adic $L$-function.  Pozzi has described explicitly the Fourier coefficients of the derivative of this family with respect to the weight variables \cite{pozzi}.  The key idea of Darmon--Pozzi--Vonk is to restrict the derivative in the anti-parallel direction along the diagonal and take the ordinary projection to obtain a classical modular form of weight $2$ for $\Gamma_0(p)$.  The idea of taking the derivative of a family of modular forms at a point of vanishing and applying a ``holomorphic projection'' operator has its roots in the seminal work of Gross--Zagier \cite{gz}, and appears more recently in Kudla's program for incoherent Eisenstein series \cite{kudla}.

Pozzi's work relates the $p$-th Fourier coefficient of this diagonal restriction to the $p$-adic logarithm of the Brumer--Stark unit $\sigma_\fb(u_p)$ for the extension $H$.  To obtain the desired weight 2 form $f_1$ on $\Gamma_0(p)$, one must take a certain linear combination with the diagonal restrictions of the two ordinary families of Eisenstein series passing through this weight 1 point.

The second form $f_2$ is defined as a generating series attached to a certain {\em rigid analytic theta cocycle}.  These are classes in $H^1(\SL_2(\Z[1/p]), \cA^*/\C_p^*)$, where $\cA^*$ denotes the group of rigid analytic nonvanishing functions on the $p$-adic upper half plane.  Darmon--Pozzi--Vonk construct classes in this space explicitly, and study their image under the {\em logarithmic annular residue map}
\[ H^1(\SL_2(\Z[1/p]), \cA^*/\C_p^*) \longrightarrow H^1(\Gamma_0(p), \Z_p). \]  They compute the spectral expansion of the form $f_2$ and thereby show that its non-constant Fourier coefficients are equal to those of $f_1$. Meanwhile, the constant coefficient is equal to the $p$-adic logarithm of $u_p(\fb)^{\an}$.  The equality of the non-constant coefficients implies that $f_1 = f_2$, and hence that the constant coefficients are equal as well, i.e.\ \[ \log_p(\sigma_\fb(u_p)) = \log_p(u_p(\fb)^{\an}) \] as desired.
It is a tantalizing problem to generalize this strategy to arbitrary totally real fields.

\section*{Acknowledgement}
We would like to thank the many mathematicians whose work has been highly influential in the development of the perspective that we have described here.  
Our work is the continuation of a long line of research connecting $L$-functions, modular forms, Galois representations, and Fitting ideals of class groups.  
In particular, we would like to thank Armand Brumer, David Burns, John Coates, Pierre Charollois, Pierre Colmez, Henri Darmon, Cornelius Greither, Benedict Gross, Masato Kurihara, Cristian Popescu, Alice Pozzi, Kenneth Ribet, Jurgen Ritter, Karl Rubin, Takamichi Sano, Michael Spiess, John Tate, Jan Vonk, Alfred Weiss, and Andrew Wiles.

\section*{Funding}
The first author is supported by a grant from the National Science Foundation, DMS-1901939.
 The second author is supported by DST-SERB grant SB/SJF/2020-21/11, SERB SUPRA grant SPR/2019/000422 and SERB MATRICS grant MTR/2020/000215.

\end{document}